\documentclass[11pt]{article}
\usepackage{amsmath}
\usepackage{amssymb}

\usepackage{theorem}
\numberwithin{equation}{section}

\newtheorem{theorem}{Theorem}[section]

\newtheorem{corollary}[theorem]{Corollary}

\newtheorem{lemma}[theorem]{Lemma}
\newtheorem{definition}{Definition}[section]

\newcommand{\ang}{{\not\negmedspace\nabla}}

\begin{document}
\title{Global well-posedness and scattering for the defocusing, energy -critical, nonlinear Schr{\"o}dinger equation in the exterior of a convex obstacle when $d = 4$}
\date{\today}
\author{Benjamin Dodson}
\maketitle

\noindent \textbf{Abstract:} In this paper we prove that the energy - critical nonlinear Schr{\"o}dinger equation in the domain exterior to a convex obstacle is globally well - posed and scattering for initial data having finite energy. To prove this we utilize frequency localized Morawetz estimates adapted to an exterior domain.

\section{Introduction}
\noindent In this paper we study the defocusing, energy - critical nonlinear Schr{\"o}dinger equation

\begin{equation}\label{1.1}
 \aligned
i u_{t} + \Delta u &= |u|^{2} u, \\
u(0,x) &= u_{0} \in \dot{H}_{0}^{1}(\Omega), \\
u|_{\partial \Omega} &= 0,
\endaligned
\end{equation}

\noindent where $\Omega = \mathbf{R}^{4} \setminus \Sigma$ is an exterior domain, $\Sigma$ is a compact, convex obstacle. We prove

\begin{theorem}\label{t1.1}
For $u_{0} \in \dot{H}_{0}^{1}(\Omega)$, $d = 4$, $(\ref{1.1})$ is globally well - posed and scattering.
\end{theorem}

\noindent As in the case when $u$ solves $i u_{t} + \Delta u = |u|^{2} u$ on $\mathbf{R}^{4}$, a solution to $(\ref{1.1})$ conserves the quantities mass,

\begin{equation}\label{1.1.1}
M(u(t)) = \int_{\Omega} |u(t,x)|^{2} dx = M(u(0)),
\end{equation}

\noindent and energy

\begin{equation}\label{1.1.2}
E(u(t)) = \frac{1}{2} \int_{\Omega} |\nabla u(t,x)|^{2} dx + \frac{1}{4} \int_{\Omega} |u(t,x)|^{4} dx = E(u(0)).
\end{equation}

\noindent $(\ref{1.1})$ is called energy critical since in $\mathbf{R}^{4}$ the symmetry $u(t,x) \mapsto \frac{1}{\lambda} u(\frac{t}{\lambda^{2}}, \frac{x}{\lambda})$ maps solutions to solutions and preserves energy. Of course, in the case of an exterior domain this scaling symmetry does not map solutions of $(\ref{1.1})$ to solutions of $(\ref{1.1})$. However, this problem behaves like the energy critical problem in $\mathbf{R}^{4}$ in many respects.\vspace{5mm}

\noindent For a domain exterior to a non - trapping obstacle \cite{IP} and \cite{BSS} proved that the quintic problem $i u_{t} + \Delta u = |u|^{4} u$ is globally well - posed and scattering for $E(u_{0})$ sufficiently small, $u$ satisfies Dirichlet or Neumann boundary conditions. \cite{LSZ} proved global well - posedness and scattering for the defocusing quintic problem when $d = 3$, $u$ is radial, and the domain is the exterior of a unit ball. The results of \cite{IP} and \cite{BSS} correspond to the results of \cite{CaWe} for the quintic problem when $d = 3$. Likewise, the techniques of \cite{LSZ} utilize the induction on energy technique used in \cite{B2} and \cite{TerryTao}.\vspace{5mm}

\begin{theorem}\label{t1.2.1}
Let $\Omega = \mathbf{R}^{d} \setminus \mathcal K$ be the exterior domain to a compact nontrapping obstacle with smooth boundary, and $\Delta$ the standard Laplace operator on $\Omega$, subject to either Dirichlet or Neumann conditions. Suppose that $p > 2$ and $q < \infty$ satisfy

\begin{equation}\label{1.1.3}
\aligned
\frac{3}{p} + \frac{n}{q} \leq \frac{n}{2}, \hspace{5mm} n = 2, \\
\frac{1}{p} + \frac{1}{q} \leq \frac{1}{2}, \hspace{5mm} n \geq 3.
\endaligned
\end{equation}

\noindent Then for the solution $v = \exp(i t \Delta) f$ to the Schr{\"o}dinger equation

\begin{equation}\label{1.1.4}
\aligned
i v_{t} + \Delta v &= 0, \\
v(0,x) &= f, \\
v|_{\partial \Omega} = 0, \hspace{5mm} &\text{ or } \hspace{5mm} \partial_{\nu} v|_{\partial \Omega} = 0,
\endaligned
\end{equation}

\noindent the following estimates hold

\begin{equation}\label{1.1.5}
\| v \|_{L_{t}^{p}([-T, T] ; L^{q}(\Omega))}	\leq C \| f \|_{\dot{H}^{s}(\Omega)},
\end{equation}

\noindent provided that

\begin{equation}\label{1.1.6}
\frac{2}{p} + \frac{n}{q} = \frac{n}{2} - s.
\end{equation}
\end{theorem}

\noindent \emph{Proof:} See \cite{BSS}. $\Box$\vspace{5mm}

\noindent In dimensions $d \geq 4$ the Strichartz estimates of \cite{BSS} are not sufficient to prove even small data global well - posedness of the energy - critical problem. Therefore, we will be content to consider the domain exterior to a convex obstacle, where we have an almost full range of Strichartz estimates.

\begin{theorem}\label{t1.3}
Suppose $u(t,x)$ is a solution to the linear Schr{\"o}dinger equation with Dirichlet boundary conditions

\begin{equation}\label{1.2}
\aligned
i u_{t} + \Delta_{D} u &= F, \\
u(0,x) &= u_{0}, \\
u|_{\partial \Omega} &= 0.
\endaligned
\end{equation}

\noindent A pair will be called admissible if $p > 2$ and

\begin{equation}\label{1.2.1}
\frac{2}{p} = d(\frac{1}{2} - \frac{1}{q}).
\end{equation}

\noindent For $(p,q)$, $(\tilde{p}, \tilde{q})$ admissible

\begin{equation}\label{1.3}
\| u \|_{L_{t}^{p}(I; L^{q}(\Omega))}	\lesssim_{p, \tilde{p}} \| u_{0} \|_{L^{2}(\Omega)}	+ \| F \|_{L_{t}^{\tilde{p}'} (I ; L^{\tilde{q}'}(\Omega))}.
\end{equation}
\end{theorem}

\noindent \emph{Proof:} See \cite{I} $\Box$.\vspace{5mm}

\noindent This theorem automatically gives small energy global well - posedness and scattering for $(\ref{1.1})$.\vspace{5mm} 

\noindent We will be able to prove theorem $\ref{t1.1}$ by utilizing the frequency truncated Morawetz estimates used in the mass - critical problem (see \cite{D2}, \cite{D3}, \cite{D5}, \cite{D4}) on $\mathbf{R}^{d}$. This technique was also used for the defocusing, energy - critical problem in $\mathbf{R}^{d}$, $d = 3, 4$. (See \cite{V2} and \cite{KV2}.) We will borrow terminology from \cite{V2} and \cite{KV2} and deal with the rapid frequency cascade and the quasi - soliton solution separately. Due to lack of scale invariance and translation invariance we will not make a concentration compactness argument. Instead, we will use induction on energy. However, the arguments used are quite reminiscent of the arguments found in \cite{D2}, \cite{D3}, \cite{D5}, \cite{D4}, \cite{V2}, and \cite{KV2}. A quick glance at \cite{V2} and \cite{KV2} will show that one might expect that the energy - critical problem in $\mathbf{R}^{4} \setminus \Sigma$ is substantially easier than the energy - critical problem in $\mathbf{R}^{3} \setminus \Sigma$. The energy - critical problem in $\mathbf{R}^{3} \setminus \Sigma$ remains out of the reach of the techniques used in this paper.\vspace{5mm}

\noindent \textbf{Function Spaces}

\noindent It will be convenient to utilize the function spaces which are a superposition of free solutions to the Schrodinger equation. See \cite{KoTa}, \cite{HHK} for more information.\vspace{5mm}

\begin{definition}\label{d1.3.1}
Let $1 \leq p < \infty$. Then $U_{\Delta_{D}}^{p}$ is an atomic space, where atoms are piecewise solutions to the linear equation $iu_{t} + \Delta_{D} u = 0$, where $\Delta_{D} = \Delta$ in the interior of $\Omega$, $\Delta_{D} = 0$ on $\partial \Omega$.

\begin{equation}\label{1.3.1}
u = \sum_{k} 1_{[t_{k}, t_{k + 1})} e^{it \Delta_{D}} u_{k}, \hspace{5mm} \sum_{k} \| u_{k} \|_{L^{2}}^{p} = 1.
\end{equation}

\noindent For any function $u$,

\begin{equation}\label{1.3.2}
\| u \|_{U_{\Delta_{D}}^{p}} = \inf \{ \sum_{\lambda} |c_{\lambda}| : u = \sum_{\lambda} c_{\lambda} u_{\lambda}, \text{$u_{\lambda}$ are $U_{\Delta_{D}}^{p}$ atoms} \}
\end{equation}

\end{definition}

\noindent For any $1 \leq p < \infty$, $U_{\Delta_{D}}^{p} \subset L^{\infty} L^{2}$. Additionally, $U_{\Delta_{D}}^{p}$ functions are continuous except at countably many points and right continuous everywhere.

\begin{definition}\label{d1.3.2}
Let $1 \leq p < \infty$. Then $V_{\Delta_{D}}^{p}$ is the space of right continuous functions $u \in L^{\infty}(L^{2})$ such that

\begin{equation}\label{1.3.3}
\| v \|_{V_{\Delta_{D}}^{p}}^{p} = \| v \|_{L^{\infty}(L^{2})}^{p} + \sup_{\{ t_{k} \} \nearrow} \sum_{k} \| e^{-it_{k} \Delta_{D}} v(t_{k}) - e^{-it_{k + 1} \Delta_{D}} v(t_{k + 1}) \|_{L^{2}}^{p}.
\end{equation}

\noindent The supremum is taken over increasing sequences $t_{k}$.
\end{definition}

\begin{theorem}\label{t1.3.3}
The function spaces $U_{\Delta_{D}}^{p}$, $V_{\Delta_{D}}^{q}$ obey the embeddings

\begin{equation}\label{1.3.4}
U_{\Delta_{D}}^{p} \subset V_{\Delta_{D}}^{p} \subset U_{\Delta_{D}}^{q} \subset L^{\infty} (L^{2}), \hspace{5mm} p < q.
\end{equation}

\noindent Let $DU_{\Delta_{D}}^{p}$ be the space of functions

\begin{equation}\label{1.3.5}
DU_{\Delta_{D}}^{p} = \{ (i \partial_{t} + \Delta_{D})u ; u \in U_{\Delta_{D}}^{p} \}.
\end{equation}

\noindent There is the easy estimate

\begin{equation}\label{1.3.6}
\| u \|_{U_{\Delta}^{p}} \lesssim \| u(0) \|_{L^{2}} + \| (i \partial_{t} + \Delta_{D}) u \|_{DU_{\Delta_{D}}^{p}}.
\end{equation}

\noindent Finally, there is the duality relation

\begin{equation}\label{1.3.7}
(DU_{\Delta_{D}}^{p})^{\ast} = V_{\Delta_{D}}^{p'}.
\end{equation}

\noindent These spaces are also closed under truncation in time.

\begin{equation}\label{2.7}
\aligned
\chi_{I} : U_{\Delta_{D}}^{p} \rightarrow U_{\Delta_{D}}^{p}, \\
\chi_{I} : V_{\Delta_{D}}^{p} \rightarrow V_{\Delta_{D}}^{p}.
\endaligned
\end{equation}
\end{theorem}

\noindent \emph{Proof:} See \cite{HHK}. $\Box$\vspace{5mm}

\noindent In particular this implies that if $(p,q)$ is an admissible pair $L_{t}^{p'} L_{x}^{q'} \subset V_{\Delta_{D}}^{2}$.\vspace{5mm}

\noindent \textbf{Remark:} From now on we will understand that $U_{\Delta}^{p}$ and $V_{\Delta}^{p}$ refers to $U_{\Delta_{D}}^{p}$ and $V_{\Delta_{D}}^{p}$ respectively.\vspace{5mm}

\noindent We have a Littlewood - Paley type theorem for an exterior domain.

\begin{theorem}\label{t1.2}
Let $\Psi \in C_{0}^{\infty}$ such that for $\lambda > 0$

\begin{equation}\label{1.1.3}
\sum_{j} \Psi(2^{-2j} \lambda) = 1.
\end{equation}

\noindent Then for $p \in (1, \infty)$, $f \in C^{\infty}(\Omega)$

\begin{equation}\label{1.1.4}
\| f \|_{L^{p}(\Omega)}	\sim_{p} \| (\sum_{j \in \mathbf{Z}} |\Psi (-2^{-2j} \Delta_{D}) f|^{2})^{1/2} \|_{L^{p}(\Omega)},
\end{equation}

\noindent and for $p \in [2, \infty)$,

\begin{equation}\label{1.1.5}
\| f \|_{L^{p}(\Omega)}	\lesssim_{p} (\sum_{j \in \mathbf{Z}} \| \Psi(-2^{2j} \Delta_{D}) f \|_{L^{p}(\Omega)}^{2})^{1/2}.
\end{equation}
\end{theorem}

\noindent \emph{Proof:} See \cite{IP2}. $\Box$\vspace{5mm}

\noindent \textbf{Remark:} As in $\mathbf{R}^{d}$ let

\begin{equation}\label{1.1.6}
u_{N} = \Psi(-N^{-2} \Delta_{D}) u,
\end{equation}

\begin{equation}\label{1.1.7}
u_{\leq 2^{j}}	= \sum_{k = -\infty}^{j} \Psi(-2^{-2k} \Delta_{D}),
\end{equation}

\begin{equation}\label{1.1.8}
u_{\leq N} + u_{> N} = u.
\end{equation}

\noindent It also follows from the fundamental theorem of calculus that

\begin{equation}\label{1.1.9}
\| u_{M} \|_{L^{\infty}(\Omega)}	\lesssim M^{d/2} \| u_{M} \|_{L^{2}(\Omega)}.
\end{equation}

\noindent As in the $\mathbf{R}^{4}$ case, to prove theorem $\ref{t1.1}$ it suffices to prove

\begin{equation}\label{1.3.1}
\| u \|_{L_{t,x}^{6}(\mathbf{R} \times \Omega)}	\leq C(E(u_{0})) < \infty.
\end{equation}

\noindent Let

\begin{equation}\label{1.9}
A(E) = \sup \{ \| u \|_{L_{t,x}^{6}(\mathbf{R} \times \Omega)}	: u \text{ solves } (\ref{1.1}), E(u(t)) = E \}.
\end{equation}

\noindent For $3 \leq d \leq 6$ it is possible to prove a stability result using exactly the same arguments which are found in \cite{TV}.

\begin{theorem}\label{t1.4}
Suppose that for $3 \leq d \leq 6$ $\tilde{u}$ is an approximate solution to $(\ref{1.1})$ in that

\begin{equation}\label{1.4}
i \tilde{u}_{t} + \Delta \tilde{u} = |\tilde{u}|^{\frac{4}{d - 2}} \tilde{u} + e,
\end{equation}

\begin{equation}\label{1.4.1}
\tilde{u}|_{\partial \Omega} = 0,
\end{equation}

\begin{equation}\label{1.5}
\| \tilde{u} \|_{L_{t,x}^{\frac{2(d + 2)}{d - 2}}(I \times \Omega)}	\leq M,
\end{equation}

\begin{equation}\label{1.6}
\| \tilde{u} \|_{L_{t}^{\infty} \dot{H}_{0}^{1}(I \times \Omega)}	\leq E,
\end{equation}

\noindent and for some $(p, q)$ admissible

\begin{equation}\label{1.7}
\| \nabla e \|_{L_{t}^{p'} L_{x}^{q'}(I \times \Omega)} \leq \epsilon,
\end{equation}

\noindent for some $\epsilon(M, E) > 0$ sufficiently small. Then there exists a solution $u(t,x)$ to $(\ref{1.1})$, $u(0,x) = \tilde{u}(0,x)$, such that for $(p, q)$ admissible

\begin{equation}\label{1.8}
\| \nabla [u - \tilde{u}] \|_{L_{t}^{p} L_{x}^{q}(I \times \Omega)}	\leq C(p, E, M) \epsilon.
\end{equation}

\end{theorem}

\noindent \emph{Proof:} We follow an argument similar to the argument in \cite{TV}. Partition $I$ into $\leq C(d,E,M)$ intervals $I_{j} = [a_{j}, b_{j}]$ such that $\| \tilde{u} \|_{L_{t,x}^{\frac{2(d + 2)}{d - 2}}(I_{j} \times \Omega)} \leq \delta$ for some small $\delta > 0$. On each $I_{j}$

\begin{equation}\label{1.9}
 \| \nabla \tilde{u} \|_{L_{t,x}^{\frac{2(d + 2)}{d}}(I_{j} \times \Omega)} \lesssim \| \tilde{u}(a_{j}) \|_{\dot{H}_{0}^{1}(\Omega)} + \| \nabla \tilde{u} \|_{L_{t,x}^{\frac{2(d + 2)}{d}}(I_{j} \times \Omega)} \| \tilde{u} \|_{L_{t,x}^{\frac{2(d + 2)}{d - 2}}(I_{j} \times \Omega)}^{\frac{4}{d - 2}} + \epsilon.
\end{equation}

\noindent This implies

\begin{equation}\label{1.10}
 \| \nabla \tilde{u} \|_{L_{t,x}^{\frac{2(d + 2)}{d - 2}}(I_{j} \times \Omega)} \lesssim E + \epsilon.
\end{equation}

\noindent Partition each $I_{j}$ into $\leq C(d,E,M)$ subintervals $I_{j,k} = [a_{j,k}, b_{j,k}]$ such that $\| \nabla \tilde{u} \|_{L_{t,x}^{\frac{2(d + 2)}{d}}(I_{j,k} \times \Omega)} \leq \delta$. Now let $v = u - \tilde{u}$. $v$ solves the initial value problem

\begin{equation}\label{1.11}
\aligned
 (i \partial_{t} + \Delta) v &= |u|^{\frac{4}{d - 2}} u - |\tilde{u}|^{\frac{4}{d - 2}} \tilde{u} - e, \\ v(0) &= 0.
\endaligned
\end{equation}

\noindent Using Strichartz estimates,

\begin{equation}\label{1.12}
\aligned
 \| \nabla v \|_{U_{\Delta}^{2}(I_{j,k} \times \Omega)} \lesssim \| v(a_{j,k}) \|_{\dot{H}_{0}^{1}(\Omega)} + \\ \| \nabla v \|_{U_{\Delta}^{2}(I_{j,k} \times \Omega)}  (\| \tilde{u} \|_{L_{t,x}^{\frac{2(d + 2)}{d - 2}}(I_{j,k} \times \Omega)} + \| \nabla \tilde{u} \|_{L_{t,x}^{\frac{2(d + 2)}{d}}(I_{j,k} \times \Omega)} + \| \nabla v \|_{L_{t,x}^{\frac{2(d + 2)}{d}}(I_{j,k} \times \Omega)})^{\frac{4}{d - 2}} + \epsilon.
\endaligned
\end{equation}

\noindent This implies that if $\| \nabla v \|_{U_{\Delta}^{2}(I_{j,k} \times \Omega)}$ is sufficiently small

\begin{equation}\label{1.13}
 \| \nabla v \|_{U_{\Delta}^{2}(I_{j,k} \times \Omega)} \lesssim \| v(a_{j,k}) \|_{\dot{H}_{0}^{1}(\Omega)} + \epsilon.
\end{equation}

\noindent For $\epsilon(E, M) > 0$ sufficiently small, since $\| v(0) \|_{\dot{H}_{0}^{1}(\Omega)} = 0$, and there are finitely many $I_{j,k}$ subintervals,

\begin{equation}\label{1.14}
 \| \nabla v \|_{U_{\Delta}^{2}(I \times \Omega)} \leq C(d,E,M) \epsilon.
\end{equation}

 $\Box$\vspace{5mm}

\noindent \textbf{Remark:} At the present time this stability result cannot be extended to $d > 6$ using the stability arguments of \cite{TV} due to the lack of exotic Strichartz estimates in a convex domain.\vspace{5mm}

\noindent By theorem $\ref{t1.4}$ $A(E)$ is a continuous function of $E$. This implies $\{ E : A(E) = \infty \}$ is a closed set, and therefore has a minimal element $E_{0}$. We prove that $E_{0} = \infty$.\vspace{5mm}

\noindent Suppose $E_{0} < \infty$. We use the bilinear virial identities of \cite{PV} to prove a bilinear Strichartz estimate for two solutions to the linear problem $i u_{t} + \Delta u = 0$, $u|_{\partial \Omega} = 0$ outside a convex obstacle. This result combined with theorem $\ref{t1.4}$ is enough to prove that a solution $u$ to $(\ref{1.1})$ with energy $E_{0}$, $\| u \|_{L_{t,x}^{6}(I \times \Omega)}	= M$, $M$ very large, must concentrate at some frequency scale $N(t)$. Partitioning $I$ into subintervals $J_{k}$ such that $\| u \|_{L_{t,x}^{6}(I \times \Omega)} = 1$, we see that $u$ must be concentrated at frequency scale $N(t) \sim N_{k}$ for some $N_{k}$. Moreover, some of the solution $u$ must be concentrated at a spatial scale $\sim \frac{1}{N_{k}}$ for length of time $\sim \frac{1}{N_{k}^{2}}$. This combined with the interaction Morawetz estimates of \cite{PV} is enough to rule out a quasi - soliton like solution. Conservation of mass rules out a rapid cascade - like solution.\vspace{5mm}

\noindent At this point it will be beneficial to say a few words about possible further developments. The purpose of this paper is two - fold. First, it is written to show that the techniques of \cite{D2}, \cite{D3}, \cite{D4}, \cite{D5}, \cite{V2}, and \cite{KV2} require very little in the way of knowledge of the fundamental solution or anything that is extremely Fourier analytic in nature.\vspace{5mm}

\noindent The second purpose is to attempt to understand the energy - critical problem in the exterior of a convex obstacle for all $d \geq 3$. The same techniques could yield global well - posedness and scattering for $d = 5$ as well. This will not be discussed in this paper because the fact that $|u|^{\frac{4}{3}} u$ is not an algebraic nonlinearity introduces some additional technical complications. The case $d = 6$ could probably be proved as well, although the proof seems to be hindered by the fact that theorem $\ref{t1.3}$ does not include endpoint Strichartz estimates. The case $d = 3$ also seems beyond the reach of the current techniques due to a heavy reliance in \cite{KV2} on Fourier - analytic techniques to obtain several key endpoint results. Extending this result to $d > 6$ would likely be far more difficult due to a lack of a stability theorem akin to theorem $\ref{t1.4}$.

\section{Morawetz Estimates}
\noindent The mainstay of the argument in this paper is the Morawetz estimates of \cite{PV} outside a star shaped obstacle. Therefore, we will summarize the argument before.

\begin{theorem}\label{t2.1}
Suppose $\Sigma$ is a compact star - shaped obstacle and $\Omega = \mathbf{R}^{d} \setminus \Sigma$ is the exterior to $\Sigma$. Let $u$ be a solution to

\begin{equation}\label{2.1}
\aligned
i u_{t} + \Delta u &= \mu |u|^{p} u, \\
u|_{t = 0} &= u_{0}.
\endaligned
\end{equation}

\noindent If $d \geq 3$, $\mu \geq 0$,

\begin{equation}\label{2.2}
\aligned
\int_{0}^{T} \int_{\partial \Omega} |\partial_{n} u|^{2} dS + \int_{0}^{T} \int_{\Omega} \frac{1}{(1 + |x|^{2})^{3/2}}	(|\nabla u(t,x)|^{2} + |u(t,x)|^{2}) dx dt	 \\+ \frac{\mu p}{p + 2} \int_{0}^{T} \int_{\Omega} \frac{1}{(1 + |x|^{2})^{1/2}} |u(t,x)|^{p + 2} dx dt	\lesssim \sup_{t \in [0, T]} \| u(t) \|_{\dot{H}_{0}^{1/2}(\Omega)}^{2}.
\endaligned
\end{equation}
\end{theorem}

\noindent \emph{Proof:} We repeat the proof found in \cite{PV}. Let $h(x) = (1 + |x|^{2})^{1/2}$,

\begin{equation}\label{2.3}
M(t) = \int_{\Omega} h(x) |u(t,x)|^{2} dx.
\end{equation}

\begin{equation}\label{2.4}
\dot{M(t)}	= 2 \int_{\Omega} h_{j}(x) Im[\bar{u} \partial_{j} u](t,x) dx.
\end{equation}

\begin{equation}\label{2.5}
\ddot{M(t)}		= \int_{\Omega} h_{j}(x) [-\partial_{k}^{2} \bar{u} \partial_{j} u + \bar{u} \partial_{j} \partial_{k}^{2} u + u \partial_{j} \partial_{k}^{2} \bar{u} - \partial_{k}^{2} u \partial_{j} \bar{u}] dx
\end{equation}

\begin{equation}\label{2.6}
- \mu \int_{\Omega} h_{j}(x) [|u|^{p} \bar{u} \partial_{j} u - \bar{u} \partial_{j} (|u|^{p} u) - u \partial_{j} (|u|^{p} \bar{u}) + |u|^{p} u \partial_{j} \bar{u}] dx.
\end{equation}

\begin{equation}\label{2.7}
\int_{\Omega} h_{j}(x) [\bar{u} \partial_{j} \partial_{k}^{2} u + u \partial_{j} \partial_{k}^{2} \bar{u}] dx
\end{equation}

\begin{equation}\label{2.8}
= -\int_{\Omega} h_{j}(x) [\partial_{k} \bar{u} \partial_{j} \partial_{k} u + \partial_{k} u \partial_{j} \partial_{k} \bar{u}] dx
\end{equation}

\begin{equation}\label{2.9}
- \int_{\Omega} h_{jk}(x) [\bar{u} \partial_{j} \partial_{k} u + u \partial_{j} \partial_{k} \bar{u}] dx.
\end{equation}

\begin{equation}\label{2.10}
(\ref{2.9}) = 2 \int_{\Omega} h_{jk}(x) Re(\partial_{k} \bar{u} \partial_{j} u)(t,x) dx	- \int_{\Omega} (\Delta \Delta h(x))	|u(t,x)|^{2} dx.
\end{equation}

\begin{equation}\label{2.10.1}
(\ref{2.8})	- \int_{\Omega} h_{j}(x) [\partial_{k}^{2} \bar{u} \partial_{j} u + \partial_{k}^{2} u \partial_{j} \bar{u}](t,x) dx	= -2 \int_{\Omega} h_{j}(x) \partial_{k} Re[(\partial_{k} \bar{u}) (\partial_{j} u)](t,x) dx
\end{equation}

\begin{equation}\label{2.10.2}
= 2 \int_{\Omega} h_{jk}(x) Re[(\partial_{k} \bar{u}) (\partial_{j} u)](t,x) dx + 2 \int_{\partial \Omega} n_{k} h_{j}(x) Re[(\partial_{k} \bar{u})(\partial_{j} u)](t,x) d\sigma(x),
\end{equation}

\noindent where $\vec{n}$ is the outward pointing unit normal to $\Sigma$, $d\sigma$ is the surface measure on $\partial \Omega$. Therefore,

\begin{equation}\label{2.11}
\aligned
(\ref{2.5})	= 4 \int_{\Omega} h_{jk}(x) Re[(\partial_{k} \bar{u})(\partial_{j} u)](t,x) dx	- \int_{\Omega} (\Delta \Delta h(x)) |u(t,x)|^{2} dx	\\ +  2 \int_{\partial \Omega} n_{k} h_{j}(x) Re[(\partial_{k} \bar{u})(\partial_{j} u)](t,x) d\sigma(x).
\endaligned
\end{equation}

\begin{equation}\label{2.12}
2 \int_{\partial \Omega} n_{k} h_{j}(x) Re[(\partial_{k} \bar{u})(\partial_{j} u)](t,x) d\sigma(x) = 2 \int_{\Omega} h_{j}(x) Re[(\partial_{n} \bar{u})(\partial_{j} u)](t,x) d\sigma(x).
\end{equation}

\noindent Because $u|_{\partial_{\Omega}} = 0$, $\nabla u = (\partial_{n} u) \vec{n}$. Therefore,

\begin{equation}\label{2.13}
(\ref{2.12})	= 2 \int_{\partial \Omega} (\partial_{n} h) |\partial_{n} u(t,x)|^{2} d\sigma(x).
\end{equation}

\noindent Let $h(x) = (1 + |x|^{2})^{1/2}$. $\nabla h(x) = \frac{x}{(1 + |x|^{2})^{1/2}}$ so $\partial_{n} h > 0$ on $\partial \Omega$ since $\Sigma$ is star - shaped.

\begin{equation}\label{2.14}
-\Delta \Delta h(x)	= (d - 1)(d - 3)(1 + |x|^{2})^{-3/2} 	+3d (1 + |x|^{2})^{-7/2}	+ 3 (d - 1)(1 + |x|^{2})^{-7/2}	+ 6(d - 3) |x|^{2} (1 + |x|^{2})^{-7/2}.
\end{equation}

\noindent $(\ref{2.14}) \geq 0$ for $d \geq 3$.

\begin{equation}\label{2.15}
4 \int_{\Omega} \partial_{k} (\frac{x_{j}}{(1 + |x|^{2})^{1/2}})	Re[(\partial_{j} \bar{u})(\partial_{k} u)](t,x) dx
\end{equation}

\begin{equation}\label{2.16}
= 4 \int_{\Omega} \frac{1}{(1 + |x|^{2})^{1/2}}	|\nabla u(t,x)|^{2} dx - 4 \int_{\Omega} \frac{|x|^{2}}{(1 + |x|^{2})^{3/2}}	|\partial_{r} u(t,x)|^{2} dx
\end{equation}

\begin{equation}\label{2.17}
= 4 \int_{\Omega} \frac{1}{(1 + |x|^{2})^{3/2}}	|\partial_{r} u(t,x)|^{2} dx	+ 4 \int_{\Omega} \frac{|x|^{2}}{(1 + |x|^{2})^{3/2}}	|\ang u(t,x)|^{2} dx.
\end{equation}

\noindent Finally, integrating by parts,

\begin{equation}\label{2.18}
(\ref{2.6})	= \frac{\mu p}{p + 2} \int ((1 + |x|^{2})^{-3/2} + (d - 1)(1 + |x|^{2})^{-1/2}) |u(t,x)|^{p + 2} dx.
\end{equation}

\noindent Combining $(\ref{2.11})$, $(\ref{2.14})$, $(\ref{2.17})$, and $(\ref{2.18})$ proves the theorem. $\Box$

\begin{theorem}\label{t2.2}
Suppose $\Omega$ is a star - shaped domain. Let $d \geq 1$, $u$, $v$ be two solutions to

\begin{equation}\label{2.19}
\aligned
i u_{t} + \Delta u &= \mu |u|^{p} u, \\
u|_{t = 0} &= u_{0},
\endaligned
\end{equation}

\begin{equation}\label{2.20}
\aligned
i v_{t} + \Delta v &= \mu |v|^{p} v, \\
v|_{t = 0} &= v_{0}.
\endaligned
\end{equation}

\noindent Let $\omega \in S^{d - 1}$.

\begin{equation}\label{2.21}
4 \int_{0}^{T} \int_{\Omega \times \Omega} |\partial_{\omega}(u(t, x_{1} + x^{\perp}) \bar{v}(t, x_{1} + y^{\perp}))|^{2} dx_{1} dx^{\perp} dy^{\perp} dt
\end{equation}

\begin{equation}\label{2.22}
+ \frac{2 \mu p}{p + 2} \int_{0}^{T} \int_{\Omega \times \Omega} |u(t, x_{1} + x^{\perp})|^{2} |v(t, x_{1} + y^{\perp})|^{p} dx_{1} dx^{\perp} dy^{\perp} dt
\end{equation}

\begin{equation}\label{2.23}
+ \frac{2 \mu p}{p + 2} \int_{0}^{T} \int_{\Omega \times \Omega} |u(t, x_{1} + x^{\perp})|^{p} |v(t, x_{1} + y^{\perp})|^{2} dx_{1} dx^{\perp} dy^{\perp} dt
\end{equation}

\begin{equation}\label{2.24}
\lesssim (\sup_{t \in [0, T]} \| u(t) \|_{\dot{H}_{0}^{1/2}(\Omega)}^{2})	(\sup_{t \in [0, T]} \| v(t) \|_{L^{2}(\Omega)}^{2})	+  (\sup_{t \in [0, T]} \| v(t) \|_{\dot{H}_{0}^{1/2}(\Omega)}^{2})	(\sup_{t \in [0, T]} \| u(t) \|_{L^{2}(\Omega)}^{2}).
\end{equation}

\end{theorem}

\noindent \textbf{Remark:} This was also proved in \cite{PV}.\vspace{5mm}

\noindent \emph{Proof:} Let

\begin{equation}\label{2.25}
I_{\omega}(\mu, u, v)	= \int_{\Omega \times \Omega} |(x - y) \cdot \omega|	|u(t,x)|^{2} |v(t,y)|^{2} dx dy.
\end{equation}

\noindent Without loss of generality suppose $\omega = (1, 0, ..., 0)$.

\begin{equation}\label{2.26}
\dot{I}_{\omega}(\mu, u, v)	= 2 \int_{\Omega \times \Omega}	\frac{(x - y)_{1}}{|(x - y)_{1}|}	Im[\bar{u} \partial_{1} u](t,x) |v(t,y)|^{2} dx dy
+ 2 \int \frac{(y - x)_{1}}{|(y - x)_{1}|}	|u(t,x)|^{2} Im[\bar{v} \partial_{1} v](t,y) dx dy.
\end{equation}

\begin{equation}\label{2.27}
\ddot{I}_{\omega}(\mu, u, v)	= \int_{\Omega \times \Omega} \frac{(x - y)_{1}}{|(x - y)_{1}|}	|v(t,y)|^{2} [-\partial_{k}^{2} \bar{u} \partial_{1} u + \bar{u} \partial_{1} \partial_{k}^{2} u + u \partial_{1} \partial_{k}^{2} \bar{u}	- \partial_{k}^{2} u \partial_{1} \bar{u}](t,x) dx dy
\end{equation}

\begin{equation}\label{2.28}
+ \int_{\Omega \times \Omega} \frac{(y - x)_{1}}{|(y - x)_{1}|}	|u(t,x)|^{2} [-\partial_{k}^{2} \bar{v} \partial_{1} v + \bar{v} \partial_{1} \partial_{k}^{2} v + v \partial_{1} \partial_{k}^{2} \bar{v}	- \partial_{k}^{2} v \partial_{1} \bar{v}](t,y) dx dy
\end{equation}

\begin{equation}\label{2.29}
- 8 \int_{\Omega \times \Omega}	Im[\bar{u} \partial_{1} u](t,x_{1} + x^{\perp})	Im[\bar{v} \partial_{1} v](t, x_{1} + y^{\perp}) dx_{1} dx^{\perp} dy^{\perp}
\end{equation}

\begin{equation}\label{2.30}
- 2 \mu \int_{\Omega \times \Omega}	\frac{(x - y)_{1}}{|(x - y)_{1}|}	[-|u|^{p} \bar{u} \partial_{1} u + \bar{u} \partial_{1}(|u|^{p} u)	+ u \partial_{1}(|u|^{p} \bar{u})	- |u|^{p} u \partial_{1} \bar{u}](t,x) |v(t,y)|^{2} dx dy
\end{equation}

\begin{equation}\label{2.31}
- 2 \mu \int_{\Omega \times \Omega}	\frac{(y - x)_{1}}{|(y - x)_{1}|}	|u(t,x)|^{2} [-|v|^{p} \bar{u} \partial_{1} v + \bar{v} \partial_{1}(|v|^{p} v)	+ v \partial_{1}(|v|^{p} \bar{v})	- |v|^{p} v \partial_{1} \bar{v}](t,y)  dx dy.
\end{equation}

\noindent \textbf{Remark:} We use the notation $x = x_{1} + x^{\perp}$, where $x_{1} = x(x \cdot (1,0,...0))$, $x^{\perp} = x - x_{1}$, and the notation

\begin{equation}\label{2.32}
\int_{\Omega \times \Omega} f(x,y) dx_{1} dx^{\perp} dy^{\perp}	= \int_{-\infty}^{\infty} \int_{x^{\perp} : x_{1} + x^{\perp} \in \Omega}	\int_{y^{\perp} : x_{1} + y^{\perp} \in \Omega}	f(x,y) dy^{\perp} dx^{\perp} dx_{1}.
\end{equation}

\noindent Following the same analysis as in the proof of theorem $\ref{t2.1}$,

\begin{equation}\label{2.33}
(\ref{2.27})	= 4	\int_{\Omega \times \Omega} \partial_{k} (\frac{(x - y)_{1}}{|(x - y)_{1}|})	Re[(\partial_{k} \bar{u})(\partial_{1} u)](t,x) |v(t,y)|^{2} dx dy
\end{equation}

\begin{equation}\label{2.34}
- \int_{\Omega \times \Omega}	\partial_{1}(\frac{(x - y)_{1}}{|(x - y)_{1}|})	\Delta(|u|^{2})(t,x) |v(t,y)|^{2} dx dy
\end{equation}

\begin{equation}\label{2.35}
+ 2 \int_{\Omega} |v(t,y)|^{2}	\int_{\partial \Omega}	\frac{(x - y)_{1}}{|(x - y)_{1}|}	Re[(\partial_{n} \bar{u})(\partial_{1} u)](t,x) dx dy.
\end{equation}

\noindent Since $\nabla u = (\partial_{n} u) \vec{n}$, by theorem $\ref{t2.1}$

\begin{equation}\label{2.36}
2 \int_{0}^{T} \int_{\Omega}	|v(t,y)|^{2}	\int_{\partial \Omega}	\frac{(x - y)_{1}}{|(x - y)_{1}|}	Re[(\partial_{n} \bar{u})(\partial_{1} u)](t,x) d\sigma(x) dy
\end{equation}

\begin{equation}\label{2.37}
\lesssim	\int_{0}^{T}	\int_{\Omega} |v(t,y)|^{2}	\int_{\partial \Omega} |\partial_{n} u(t,x)|^{2} d\sigma(x) dy	\lesssim	(\sup_{t \in [0, T]} \| u(t) \|_{\dot{H}_{0}^{1/2}(\Omega)}^{2})	(\sup_{t \in [0, T]} \| v(t) \|_{L^{2}(\Omega)}^{2}).
\end{equation}

\noindent This takes care of $(\ref{2.35})$. Next,

\begin{equation}\label{2.38}
\aligned
(\ref{2.34})	=	-\int_{\Omega \times \Omega}	\Delta_{x}(|u(t,x_{1} + x^{\perp})|^{2})	|v(t,x_{1} + y^{\perp})|^{2} dx_{1} dx^{\perp}	dy^{\perp}		\\= \int_{\Omega \times \Omega}	\partial_{x_{1}}(|u(t,x_{1} + x^{\perp})|^{2})	\partial_{x_{1}}(|v(t, x_{1} + y^{\perp})|^{2})	dx_{1} dx^{\perp} dy^{\perp}.
\endaligned
\end{equation}

\noindent Finally,

\begin{equation}\label{2.39}
(\ref{2.33})	= 4 \int_{\Omega \times \Omega}	|\partial_{1} u(t,x + x^{\perp})|^{2} |v(t, x_{1} + y^{\perp})|^{2} dx_{1} dx^{\perp} dy^{\perp}.
\end{equation}

\noindent Make an identical argument for $(\ref{2.28})$.

\begin{equation}\label{2.40}
4 |\partial_{1} u|^{2} |v|^{2}	+ 4 |\partial_{1} v|^{2} |u|^{2}	- 8 Im[\bar{u} \partial_{1} u] Im[\bar{v} \partial v]		+ 2 \partial_{1}(|u|^{2}) \partial_{1}(|v|^{2})		= 4 |\partial_{1}(u \bar{v})|^{2}.
\end{equation}

\noindent By the fundamental theorem of calculus,

\begin{equation}\label{2.41}
\sup_{t \in [0, T]} \dot{I}_{\omega}(\mu, u, v)	\lesssim  (\sup_{t \in [0, T]} \| u(t) \|_{\dot{H}_{0}^{1/2}(\Omega)}^{2})	(\sup_{t \in [0, T]} \| v(t) \|_{L^{2}(\Omega)}^{2})	+  (\sup_{t \in [0, T]} \| v(t) \|_{\dot{H}_{0}^{1/2}(\Omega)}^{2})	(\sup_{t \in [0, T]} \| u(t) \|_{L^{2}(\Omega)}^{2}),
\end{equation}

\noindent $(\ref{2.37})$, and integrating $(\ref{2.30})$ and $(\ref{2.31})$ by parts, the proof is complete. $\Box$

\begin{corollary}\label{c2.3}
When $d = 1$,

\begin{equation}\label{2.42}
\int_{0}^{T}	\int_{\Omega} |\partial_{x}(u \bar{v})|^{2}(t,x) + \frac{2 \mu p}{p + 2} |u|^{2} |v|^{p + 2}(t,x)	+ \frac{2 \mu p}{p + 2} |u|^{p + 2} |v|^{2}(t,x) dx dt
\end{equation}

\begin{equation}\label{2.43}
\lesssim	(\sup_{t \in [0, T]} \| u(t) \|_{\dot{H}_{0}^{1/2}(\Omega)}^{2})	(\sup_{t \in [0, T]} \| v(t) \|_{L^{2}(\Omega)}^{2})	+ (\sup_{t \in [0, T]} \| v(t) \|_{\dot{H}_{0}^{1/2}(\Omega)}^{2})	(\sup_{t \in [0, T]} \| u(t) \|_{L^{2}(\Omega)}^{2}).
\end{equation}

\end{corollary}

\noindent Now replace $|(x - y) \cdot \omega$ with a more general $\rho(x - y)$ with positive definite Hessian $H_{\rho}(x - y)$.

\begin{theorem}\label{t2.5}
Let $F(u,v)(x,y) = \bar{v}(y) \nabla_{x} u(x) + u(x) \nabla_{y} \bar{v}(y)$ and $G(u,v)(x,y) = v(y) \nabla_{x} u(x) - u(x) \nabla_{y} v(y)$. Then

\begin{equation}\label{2.46}
4 \int_{0}^{T} \int_{\Omega \times \Omega}	H_{\rho}(x - y) (F(u,v)(x,y), \bar{F}(u,v)(x,y)) dx dy dt
\end{equation}

\begin{equation}\label{2.47}
+ \frac{\mu p}{p + 2} \int_{0}^{T}	\int_{\Omega \times \Omega}	|v|^{2}(t,y) (\Delta_{x} \rho)(x - y) |u|^{p + 2}(t,x) dx dy dt
\end{equation}

\begin{equation}\label{2.48}
+ \frac{\mu p}{p + 2} \int_{0}^{T}	\int_{\Omega \times \Omega}	|v|^{p + 2}(t,y) (\Delta_{x} \rho)(x - y) |u|^{2}(t,x) dx dy dt
\end{equation}

\begin{equation}\label{2.49}
\lesssim	(\sup_{t \in [0, T]} \| u(t) \|_{\dot{H}_{0}^{1/2}(\Omega)}^{2})	(\sup_{t \in [0, T]} \| v(t) \|_{L^{2}(\Omega)}^{2})	+  (\sup_{t \in [0, T]} \| v(t) \|_{\dot{H}_{0}^{1/2}(\Omega)}^{2})	(\sup_{t \in [0, T]} \| u(t) \|_{L^{2}(\Omega)}^{2}).
\end{equation}

\end{theorem}

\noindent \textbf{Remark:} This also appears in \cite{PV}.\vspace{5mm}

\noindent \emph{Proof:} Follow the proof of theorem $\ref{t2.2}$.

\begin{equation}\label{2.50}
\aligned
4 \rho_{jk}(x - y) Re((\partial_{j} \bar{u})(\partial_{k} u))(t,x)	|v(t,y)|^{2}		+ 4 \rho_{jk}(x - y) Re((\partial_{j} \bar{v})(\partial_{k} v))(t,y) |u(t,x)|^{2}	\\
- 8 \rho_{jk}(x - y) Im[\bar{u} \partial_{j} u](t,x) Im[\bar{v} \partial_{k} v](t,y) - 2 \rho_{jk}(x - y) \partial_{j}(|u|^{2})(t,x) \partial_{k}(|v|^{2})(t,y)	\\
= 4 H_{\rho}(x - y)(F(u,v), \bar{F}(u,v))	\\= 4 H_{\rho}(x - y)(G(u,v), \bar{G}(u,v))	+ \Delta \rho(x - y) \nabla_{x}(|u|^{2})(t,x) \cdot \nabla_{y}(|v|^{2})(t,y).
\endaligned
\end{equation}

\noindent We use these arguments to prove a bilinear Strichartz estimate in the exterior of a star - shaped domain.

\begin{theorem}\label{t2.6}
Suppose $\Omega$ is a domain exterior to a star shaped obstacle. Suppose $u_{0} = \Psi(-M^{-2} \Delta_{D}) u_{0}$ and $v_{0} = \Psi(-N^{-2} \Delta_{D}) v_{0}$, $M \leq N$. If $u$ and $v$ are linear solutions to

\begin{equation}\label{2.51}
\aligned
i u_{t} + \Delta u &= 0, \\
u(0) &= u_{0},
\endaligned
\end{equation}

\begin{equation}\label{2.52}
\aligned
i v_{t} + \Delta v &= 0, \\
v(0) &= v_{0},
\endaligned
\end{equation}

\noindent then

\begin{equation}\label{2.52.1}
\| \nabla (u \bar{v})	\|_{L_{t,x}^{2}(\mathbf{R} \times \Omega)}		\lesssim	M^{(d - 1)/2} N^{1/2}	\| u_{0} \|_{L^{2}(\Omega)} \| v_{0} \|_{L^{2}(\Omega)}.
\end{equation}
\end{theorem}

\noindent \emph{Proof:} By elementary Strichartz estimates the theorem follows for $M \sim N$. By the fundamental theorem of calculus, when $d = 1$, $\| u \|_{L^{\infty}(\mathbf{R})}^{2}	\lesssim	\| \nabla u \|_{L^{2}(\mathbf{R})} \| u \|_{L^{2}(\mathbf{R})}$. In general,

\begin{equation}\label{2.53}
\| u \|_{L^{\infty}(\Omega)}^{2}	\lesssim	\| (-\Delta_{D})^{d/2} u \|_{L^{2}(\Omega)} \| u \|_{L^{2}(\Omega)}	\lesssim	\frac{1}{M^{d}}	\| (-\Delta_{D})^{d/2} u \|_{L^{2}(\Omega)}^{2} + M^{d} \| u \|_{L^{2}(\Omega)}^{2}.
\end{equation}

\noindent Let $\chi \in C_{0}^{\infty}(\mathbf{R}^{d})$, $\chi = 1$ on $|x| \leq \frac{1}{2}$, $\chi = 0$ on $|x| > 1$. For $x_{0} \in \Omega$,

\begin{equation}\label{2.54}
|u(x_{0})|^{2}	\lesssim \frac{1}{M^{d}}	\| \chi(\frac{M(x - x_{0})}{2}) ((-\Delta_{D})^{d/2} u) \|_{L^{2}(\Omega)}^{2}	+ M^{d} \| \chi(\frac{M(x - x_{0})}{2}) u \|_{L^{2}(\Omega)}^{2}.
\end{equation}

\noindent Making basic Strichartz estimates,

\begin{equation}\label{2.55}
M^{d} \int_{0}^{T} \int_{\Omega} |u(t, x_{1} + x^{\perp})|^{2}	\int_{|\tau| \leq \frac{2}{M}}		\int_{|y^{\perp} - x^{\perp}| \leq \frac{2}{M}}	|\partial_{1} v(t, x_{1} + \tau e_{1} + y^{\perp})|^{2} dy^{\perp} d\tau dx^{\perp} dx_{1}
\end{equation}

\begin{equation}\label{2.56}
\lesssim M^{d - 1} \| u_{0} \|_{L^{2}(\Omega)}^{2} \| v_{0} \|_{L^{2}(\Omega)}^{2}.
\end{equation}

\noindent Therefore, by theorem $\ref{t2.2}$, $(\ref{2.56})$, combined with the fact that $|(x - y)_{1} - \tau|$ has a positive definite Hessian,

\begin{equation}\label{2.57}
M^{d}	\int_{0}^{T} \int_{\Omega} |\partial_{1} (u(t, x_{1} + x^{\perp}))|^{2} \int_{|\tau| \leq \frac{2}{M}} \int_{|y^{\perp} - x^{\perp}| \leq \frac{2}{M}}	|v(t,x_{1} + \tau e_{1} + y^{\perp})|^{2} dy^{\perp} d\tau dx^{\perp} dx_{1} dt
\end{equation}

\begin{equation}\label{2.58}
\aligned
\lesssim	M^{d - 1} (\sup_{t \in [0, T]} \| u(t) \|_{L^{2}(\Omega)}^{2}) (\sup_{t \in [0, T]} \| v(t) \|_{\dot{H}_{0}^{1/2}(\Omega)}^{2})	\\+ M^{d - 1} (\sup_{t \in [0, T]} \| v(t) \|_{L^{2}(\Omega)}^{2}) (\sup_{t \in [0, T]} \| u(t) \|_{\dot{H}_{0}^{1/2}(\Omega)}^{2})	\lesssim	M^{d - 1} N \| u_{0} \|_{L^{2}(\Omega)}^{2} \| v_{0} \|_{L^{2}(\Omega)}^{2}.
\endaligned
\end{equation}

\noindent Similarly, since $\Delta_{D}$ commutes with the solution operator to $(\ref{2.52})$,

\begin{equation}\label{2.59}
\frac{1}{M^{d}}	\int_{0}^{T} \int_{\Omega}	|\partial_{1} (u(t, x_{1} + x^{\perp}))|^{2}	\int_{|\tau| \leq \frac{2}{M}}	\int_{|y^{\perp} - x^{\perp}| \leq \frac{2}{M}}	|(-\Delta_{D})^{d/2} v(t, x_{1} + \tau e_{1} + y^{\perp})|^{2} dy^{\perp} d\tau dx^{\perp} dx_{1} dt
\end{equation}

\begin{equation}\label{2.60}
\aligned
\lesssim	\frac{1}{M^{d + 1}}	(\sup_{t \in [0, T]} \| u(t) \|_{L^{2}(\Omega)}^{2}) (\sup_{t \in [0, T]} \| (-\Delta_{D})^{d/2}	v(t) \|_{\dot{H}_{0}^{1/2}(\Omega)}^{2})	\\+ \frac{1}{M^{d + 1}} (\sup_{t \in [0, T]} \|(-\Delta_{D})^{d/2} v(t) \|_{L^{2}(\Omega)}^{2}) (\sup_{t \in [0, T]} \| u(t) \|_{\dot{H}_{0}^{1/2}(\Omega)}^{2})	\lesssim	M^{d - 1} N \| u_{0} \|_{L^{2}(\Omega)}^{2} \| v_{0} \|_{L^{2}(\Omega)}^{2}.
\endaligned
\end{equation}

\noindent This completes the proof of theorem $\ref{t2.6}$. $\Box$\vspace{5mm}

\noindent We can now prove that a solution to $(\ref{1.1})$ with energy $E_{0}$ and very high $\| u \|_{L_{t,x}^{6}(J \times \Omega)}$ norm for some compact interval $J$ must concentrate in frequency.\vspace{5mm}

\noindent By theorem $\ref{t2.6}$ and theorem $\ref{t1.4}$,

\begin{theorem}\label{t2.7}
 Suppose $u$ solves $(\ref{1.1})$, $E(u(t)) = E_{0}$,

\begin{equation}\label{2.61}
\| u \|_{L_{t,x}^{6}(J \times \Omega)}	= M,
\end{equation}

\noindent for some $M$ very large. Then fix $\delta(E_{0}) > 0$ sufficiently small and partition $J$ into subintervals $J_{k}$ such that

\begin{equation}\label{2.62}
 \| u \|_{L_{t,x}^{6}(J_{k} \times \Omega)}	= \delta.
\end{equation}

\noindent For each $J_{k}$ there exists $N_{k} \in (0, \infty)$ such that for all $\eta > \eta_{0}(M)$, $\eta_{0}(M) \searrow 0$ as $M \nearrow \infty$, there exists $C(\eta) < \infty$ such that

\begin{equation}\label{2.63}
\| P_{\leq \frac{1}{C(\eta)} N_{k}}	u \|_{L_{t}^{\infty} \dot{H}^{1}(\Omega)(J_{k} \times \Omega)}	+ \| P_{\leq \frac{1}{C(\eta)} N_{k}}	u \|_{L_{t}^{\infty} L_{x}^{4}(J_{k} \times \Omega)} < \eta.
\end{equation}

\begin{equation}\label{2.64}
 \| P_{\geq C(\eta) N_{k}} u \|_{L_{t}^{\infty} \dot{H}^{1}(J_{k} \times \Omega)} + \| P_{\geq C(\eta) N_{k}}	u \|_{L_{t}^{\infty} L_{x}^{4}(J_{k} \times \Omega)} < \eta.
\end{equation}

\noindent Moreover, if $J_{k}$ and $J_{k + 1}$ are adjacent intervals then $N_{k} \sim_{E_{0}} N_{k + 1}$.

\end{theorem}
 
 \noindent We divide into two cases. Let $C_{0} = \inf_{t \in J} N(t)$, $M$ very large. We treat the case
 
 \begin{equation}\label{2.65}
C_{0}^{3} \sum_{J_{k} \subset J}	\frac{1}{N_{k}^{3}}	\leq K_{0}
\end{equation}

\noindent as \cite{V2}, \cite{KV2} treated the rapid cascade. We treat

\begin{equation}\label{2.66}
C_{0}^{3} \sum_{J_{k} \subset J}	\frac{1}{N_{k}^{3}}	\geq K_{0}
\end{equation}
 
 \noindent as \cite{V2}, \cite{KV2} treated the pseudo - soliton, for some $K_{0}$ to be specified later. We prove,
 
 \begin{theorem}\label{t2.8}
 There does not exist a solution to $(\ref{1.1})$,
 
 \begin{equation}\label{2.67}
 \| u \|_{L_{t,x}^{6}(J \times \Omega)}	= M,
 \end{equation}
 
 \noindent $E(u(t)) = E_{0}$, for $M$ very large.
 \end{theorem}

 \section{Long time Strichartz estimates for the rapid cascade}
 
We first rule out a scenario similar to the rapid frequency cascade. Fix $K$.

\begin{theorem}\label{t3.0}
There is a constant $K_{0}(M)$ such that $K_{0}(M) \nearrow \infty$ as $M \rightarrow \infty$ and there does not exist a solution to $(\ref{1.1})$ with energy $E(u(t)) = E_{0}$, an interval $J$ with

\begin{equation}\label{3.0}
 C_{0}^{3}		\sum \frac{1}{N_{k}^{3}} \leq K_{0}(M),
\end{equation}

\begin{equation}\label{3.0.1}
\| u \|_{L_{t,x}^{6}(J \times \Omega)}	= M.
\end{equation}
\end{theorem}

\noindent To prove this theorem we utilize a slight modification of the arguments of \cite{V2}, \cite{KV2} for the energy - critical problem in $\mathbf{R}^{d}$. See also \cite{D2}, \cite{D3}, and \cite{D5} for induction on frequency in the mass - critical case.\vspace{5mm}

\noindent As in the case of the energy - critical nonlinear Schr{\"o}dinger equation on flat space we will rule out a sufficiently large blowup solution.

\begin{theorem}\label{t3.1}
 Suppose $J$ is a union of subintervals $J_{k}$ such that for some $\epsilon > 0$

\begin{equation}\label{3.1}
 \sum \frac{1}{N_{k}^{2(2 - \epsilon)}}	= K.
\end{equation}

\noindent Then if $(p,q)$ satisfies $(\ref{1.2.1})$,

\begin{equation}\label{3.2}
 \| \nabla P_{\leq N} u \|_{L_{t}^{p} L_{x}^{q}(J \times \Omega)}	\lesssim_{\epsilon, p,q,d}	1 + K^{1/p} N^{\frac{2(2 - \epsilon)}{p}}.
\end{equation}

\end{theorem}

\noindent \emph{Proof:} Fix $\epsilon > 0$. For each dyadic $N$ we partition $J$ at level $N$. We call these corresponding intervals $J_{N}^{l}$. If $N_{k} > \frac{c}{N}$, for some small, fixed $c > 0$ to be specified later, then we say $J_{k}$ is a bad interval, $J_{N,b}^{l}$. We group the remaining $J_{k}$ subintervals into good $J_{N,g}^{l}$ intervals such that each good interval satisfies

\begin{equation}\label{3.3}
\frac{1}{2} < N^{2(2 - \epsilon)}	\sum_{J_{k} \subset J_{N,g}^{l}}	\frac{1}{N_{k}^{2(2 - \epsilon)}}	\leq 1,	
\end{equation}

\noindent or 

\begin{equation}\label{3.4}
 N^{2(2 - \epsilon)}	\sum_{J_{k} \subset J_{N,g}^{l}}	\frac{1}{N_{k}^{2(2 - \epsilon)}}	\leq \frac{1}{2},
\end{equation}

\noindent and $J_{N,g}^{l}$ is adjacent to a bad interval. It suffices to prove

\begin{lemma}\label{l3.2}
For any dyadic integer $N$ and any interval $J_{N}^{l}$,

\begin{equation}\label{3.5}
 \| \nabla P_{\leq N} u \|_{U_{\Delta}^{2}(J_{N}^{l} \times \Omega)}	\lesssim_{E_{0}, \epsilon} 1.
\end{equation}

\end{lemma}

\noindent Indeed, for any $p > 2$,

\begin{equation}\label{3.6}
 \| \nabla P_{\leq N} u \|_{L_{t}^{p} L_{x}^{q}(J_{N}^{l} \times \Omega)}	\lesssim	\| \nabla P_{\leq N} u \|_{U_{\Delta}^{2}(J_{N}^{l} \times \Omega)}	\lesssim_{E_{0}, \epsilon} 1.
\end{equation}

\noindent Since $\sharp \{ J_{N}^{l} \} \lesssim N^{2(2 - \epsilon)} K$, summing up the norm of $J_{N}^{l}$ intervals in $l^{p}$ proves theorem $\ref{t3.1}$. $\Box$\vspace{5mm}

\noindent \emph{Proof of lemma $\ref{l3.2}$:} This follows from Duhamel's formula. For $t_{0,N}^{l} \in J_{N}^{l}$, a solution to $(\ref{1.1})$ satisfies

\begin{equation}\label{3.7}
 e^{i(t - t_{0, N}^{l}) \Delta} u(t_{0,N}^{l})	+ \int_{t_{0,N}^{l}}^{t} e^{i(t - \tau) \Delta} F(u(\tau)) d\tau,
\end{equation}

\noindent where $F(u) = |u|^{2} u$. By conservation of energy, $\| \nabla e^{i(t - t_{0,N}^{l}) \Delta} u(t_{0,N}^{l}) \|_{U_{\Delta}^{2}(J_{N}^{l} \times \Omega)}	\lesssim 1$.\vspace{5mm}

\noindent Next, for a bad interval $J_{N,b}^{l}$,

\begin{equation}\label{3.8}
 \| \nabla F(u) \|_{L_{t}^{p'} L_{x}^{q'}(J_{N,b}^{l} \times \Omega)}	\lesssim	\| \nabla F(u) \|_{L_{t,x}^{3}(J_{N,b}^{l} \times \Omega)}	\| u \|_{L_{t,x}^{6}(J_{N,b}^{l} \times \Omega)}^{2}	\lesssim 1.
\end{equation}

\noindent For $N = N_{max} = \sup_{k} N_{k}$ all intervals are bad, and so we are done. We therefore proceed by induction. Let $C$ be some large, fixed constant, $c(C) > 0$ a small constant to be chosen momentarily.

\begin{equation}\label{3.9}
 \| \nabla F(P_{< CN} u) \|_{L_{t}^{2-} L_{x}^{4/3+}(J_{N,g}^{l} \times \Omega)}	\lesssim	
\end{equation}

\begin{equation}\label{3.10}
\| P_{\leq CN} u \|_{L_{t}^{2+} L_{x}^{4-}(J_{N, g}^{l} \times \Omega)}^{1+} \| P_{< cN_{k}} u \|_{L_{t}^{\infty} L_{x}^{4}(J_{N, g}^{l} \times \Omega)}^{2-} \lesssim (C)^{1+} \eta^{2-}.
\end{equation}

\noindent Next let $p = 2+$, $(p, q)$ satisfies $(\ref{1.2.1})$, $p \searrow 2$ as $\epsilon \searrow 0$. By induction, for $M > CN$, interpolating $U_{\Delta}^{2}$ estimates with $L_{t}^{\infty} \dot{H}^{1}$ estimates,

\begin{equation}\label{3.11}
\| u_{M} \|_{L_{t}^{p} L_{x}^{q}(J_{N,g}^{l} \times \Omega)}^{1+} \| u_{M} \|_{L_{t}^{\infty} L_{x}^{2}(J_{N,g}^{l} \times \Omega)}^{1-} \| u_{M} \|_{L_{t}^{\infty} L_{x}^{4}(J_{N,g}^{l} \times \Omega)}	\lesssim	\frac{1}{M^{2}}	(\frac{M}{N})^{(2 - \epsilon)(1+)} (\eta(c))^{2-}.
\end{equation}

\noindent Therefore, by Sobolev embedding,

\begin{equation}\label{3.12}
N \| P_{< N} F( u_{CN < \cdot < cN_{k}}) \|_{L_{t}^{p'} L_{x}^{q'}(J_{N,g}^{l} \times \Omega)}	\lesssim	(\eta(c))^{2-}.
\end{equation}

\noindent Choosing $\eta$ sufficiently small, $c(\eta) > 0$ closes the induction. Next, for some $\delta(p, \epsilon) > 0$,

\begin{equation}\label{3.13}
\aligned
N^{2}	\| u_{> cN_{k}} \|_{L_{t}^{p} L_{x}^{q}(J_{N,g}^{l} \times \Omega)}^{1 +\delta} \| u_{> cN_{k}} \|_{L_{t}^{\infty} L_{x}^{2}(J_{N,g}^{l} \times \Omega)}^{1 - \delta} \| u_{> cN_{k}} \|_{L_{t}^{\infty} L_{x}^{4}(J_{N,g}^{l} \times \Omega)}	\\
	\lesssim	N^{2}	(\sum \frac{1}{(cN_{k})^{2(2 - \frac{\epsilon}{2})}})^{\frac{2}{2 - \frac{\epsilon}{2}}}	\lesssim \frac{1}{c^{2}}.
\endaligned
\end{equation}

\noindent Notice that this term does not depend on the inductive hypotheses. The last inequality follows from the fact that $\frac{N}{N_{k}} \leq 1$ on good intervals. Finally, by the inductive hypothesis

\begin{equation}\label{3.14}
N \| |u_{> CN} | |u_{< CN}|^{2}	+ |u_{> CN}|^{2} |u_{< CN}|	\|_{L_{t}^{\tilde{p}'} L_{x}^{\tilde{q}'}(J_{N, g}^{l} \time \Omega)}
\end{equation}

\begin{equation}\label{3.14.1}
\aligned
\lesssim	N	\| u_{< CN} \|_{L_{t}^{2+} L_{x}^{\infty-}(J_{N,g}^{l} \times \Omega)}^{1+} \| u_{< CN} \|_{L_{t}^{\infty} L_{x}^{4}(J_{N,g}^{l} \times \Omega)}^{1-} \| u_{> CN} \|_{L_{t}^{\infty} L_{x}^{2}(J_{N,g}^{l} \times \Omega)}	\\
+ N^{2}	\| u_{> CN} \|_{L_{t}^{2+} L_{x}^{4-}(J_{N,g}^{l} \times \Omega)}^{1+}	\| u_{> CN} \|_{L_{t}^{\infty} L_{x}^{2}(J_{N,g}^{l} \times \Omega)}^{1-}	\| u_{< CN} \|_{L_{t}^{\infty} L_{x}^{4}(J_{N,g}^{l} \times \Omega)}	\lesssim	C^{1+} \eta^{1-}	+ \eta \frac{1}{c^{2}}.
\endaligned
\end{equation}

\noindent Therefore we have proved

\begin{equation}\label{3.15}
\| \nabla P_{< N} [F(u) - F(u_{< CN})] \|_{L_{t}^{p'} L_{x}^{q'}(J_{N, g}^{l} \times \Omega)}	\lesssim 1.
\end{equation}

\noindent This completes the proof of lemma $\ref{l3.2}$. $\Box$\vspace{5mm}

\noindent \textbf{Remark:} As in \cite{D2} we can upgrade $(\ref{3.2})$ to

\begin{equation}\label{3.20}
 \| \nabla P_{\leq N} u \|_{L_{t}^{p} L_{x}^{q}(J \times \Omega)}	\lesssim_{\epsilon, p,q,d}	\sigma_{J}(N) + K^{1/p} N^{\frac{2(1 + \frac{4}{d - 2} - \epsilon)}{p}},
\end{equation}

\noindent where $\sigma_{J}(N)$ is a frequency envelope that majorizes

\begin{equation}\label{3.21}
 \inf_{t \in J} \| \nabla P_{\leq N} u(t) \|_{L_{t}^{\infty} L_{x}^{2}(J \times \Omega)}.
\end{equation}

\noindent Let

\begin{equation}\label{3.22}
\sigma_{J}(N) = \sum_{j = -\infty}^{\infty} 2^{-|j|} \inf_{t \in J} \| P_{2^{j} N} u(t) \|_{L_{t}^{\infty} \dot{H}_{x}^{1}(J \times \Omega)}.
\end{equation}

\noindent This is enough to prove theorem $\ref{t3.0}$.\vspace{5mm}

\noindent \emph{Proof of theorem $\ref{t3.0}$:}	 Let $u$ be a solution to $(\ref{1.1})$, $\| u \|_{L_{t,x}^{6}(J \times \Omega)} = M$, $K = \frac{K_{0}}{C_{0}^{3}}$. For any $N$ let $u_{l} = P_{\leq N} u$, $u_{l} + u_{h} = u$. If $K N^{3} \leq 1$, then as in the previous section

\begin{equation}\label{3.23}
 \frac{d}{dt} \langle \nabla P_{\leq N} u, \nabla P_{\leq N} u \rangle	= 2i \langle \nabla P_{\leq N} F(u), \nabla P_{\leq N} u \rangle	\lesssim \sigma_{J}(N)^{2} + K N^{3}.
\end{equation}

\noindent This implies that for any $t \in J$,

\begin{equation}\label{3.24}
 u(t) = v(t) + w(t),
\end{equation}

\noindent where 

\begin{equation}\label{3.25}
\| v(t) \|_{L_{t}^{\infty} \dot{H}_{x}^{-1/4}(J \times \Omega)}	\lesssim \frac{K_{0}^{5/12}}{C_{0}^{5/4}},
\end{equation}

\noindent and $\| w(t) \|_{L_{t}^{\infty} \dot{H}_{x}^{1}(J \times \Omega)}	\searrow 0$ as $M \nearrow \infty$. However, by interpolation this implies that for $t \in J_{k}$, $\eta > \eta_{0}(M)$,

\begin{equation}\label{3.26}
\| v(t) \|_{L_{x}^{2}(\Omega)}	\lesssim \frac{C(\eta)}{N_{k}}	+ \eta^{1/5} \frac{K_{0}^{1/3}}{C_{0}}.
\end{equation}

\noindent Since $J$ satisfies $(\ref{3.0})$, this implies that for $M$ sufficiently large there exists an interval $J_{k}$ such that for $t \in J_{k}$, for any $\eta > \eta(M) > 0$,

\begin{equation}\label{3.27}
 \| v(t) \|_{L_{x}^{2}(\Omega)}	<< \frac{1}{C_{0}}.
 \end{equation}

\noindent Moreover, for some $\tilde{t}$ in this interval, with $t_{0}$ satisfying $N(t_{0}) = \inf_{t \in J} N(t)$, and without loss of generality $t_{0} < \tilde{t}$,

\begin{equation}\label{3.28}
\| u \|_{L_{t,x}^{6}([t_{0}, \tilde{t}] \times \Omega)}	\lesssim_{K_{0}} 1.
\end{equation}

\noindent By the perturbation lemma this contradicts theorem $\ref{t2.7}$. $\Box$\vspace{5mm}

\section{Interaction Morawetz estimates}

\begin{theorem}\label{t4.0}
There is a fixed constant $K_{0} < \infty$ such that for $M$ sufficiently large, there does not exist a solution to $(\ref{1.1})$ satisfying $E(u(t)) = E_{0}$,

\begin{equation}\label{4.1}
C_{0}^{3} \sum_{J_{k} \subset J}	\frac{1}{N_{k}^{3}}	\geq K_{0},
\end{equation}

\noindent and

\begin{equation}\label{4.2}
\| u \|_{L_{t,x}^{6}(J \times \Omega)}	= M.
\end{equation}
\end{theorem}

\noindent Recall that $K = \frac{K_{0}}{C_{0}^{3}}$.

\begin{theorem}\label{t4.1}
\noindent Let $u_{h} = P_{\geq K^{-1/3}} u$, $u_{h} + u_{l} = u$.

\begin{equation}\label{4.3}
\int_{J} \int_{\partial \Omega}	|\partial_{n} u_{h}(t,x)|^{2}		dS_{x} dt	\lesssim K^{1/3}.
\end{equation}
\end{theorem}

\noindent \emph{Proof:} By $\S 3$, for $(p, q)$ satisfying $(\ref{1.2.1})$, $p > 2$,

\begin{equation}\label{4.4}
\| \nabla u_{l}	\|_{L_{t}^{p} L_{x}^{q}(J \times \Omega)}	\lesssim 1.
\end{equation}

\noindent We will also postpone the proof of the estimate

\begin{equation}\label{4.5}
\| u_{h} \|_{L_{t,x}^{3}(J \times \Omega)}	+ \| u_{h} \|_{L_{t}^{\infty} L_{x}^{2}(J \times \Omega)}	\lesssim K^{1/3}.
\end{equation}

\noindent Let

\begin{equation}\label{4.6}
M(t)	= \int_{\Omega}	\frac{x_{j}}{(1 + |x|^{2})^{1/2}}	Im[\bar{u}_{h} \partial_{j} u_{h}](t,x) dx.
\end{equation}

\begin{equation}\label{4.7}
\sup_{t \in J}	|M(t)|		\lesssim K^{1/3}.
\end{equation}

\noindent Therefore, following the analysis in \cite{PV},

\begin{equation}\label{4.8}
\int_{J} \int_{\partial \Omega}	|\partial_{n} u_{h}(t,x)|^{2} dS_{x} dt	+ \int_{J} \int	\frac{1}{(1 + |x|^{2})^{3/2}}	(|\nabla u_{h}(t,x)|^{2}	+ |u_{h}(t,x)|^{2}) dx dt
\end{equation}

\begin{equation}\label{4.9}
+ \int_{J} \int_{\Omega}	\frac{x_{j}}{(1 + |x|^{2})^{1/2}}	\{ P_{h} F(u), u_{h} \}_{j} dx dt	\lesssim K^{1/3},
\end{equation}

\noindent where

\begin{equation}\label{4.10}
\{ u, v \}_{j} = Re [\bar{u} \partial_{j} v - \bar{v} \partial_{j} u].
\end{equation}

\begin{equation}\label{4.11}
\{ P_{h} F(u), u_{h} \}_{j}	= \{ F(u), u \}_{j}	- \{ F(u_{l}), u_{l} \}_{j}	- \{ F(u) - F(u_{l}), u_{l} \}_{j}	- \{ P_{l} F(u), u_{h} \}_{j}.
\end{equation}

\begin{equation}\label{4.11.1}
\{ F(u), u \}_{j}	- \{ F(u_{l}), u_{l} \}_{j}	= \frac{1}{2} \partial_{j} [|u(t,x)|^{4} - |u_{l}(t,x)|^{4}].
\end{equation}

\begin{equation}\label{4.12}
\int_{J} \int_{\Omega}	\frac{x_{j}}{(1 + |x|^{2})^{1/2}}	[F(u) - F(u_{l})] (\partial_{j} u_{l})(t,x) dx dt
\end{equation}

\begin{equation}\label{4.13}
\lesssim	\| \nabla u_{l} \|_{L_{t}^{3} L_{x}^{12}}	\| u_{h} \|_{L_{t,x}^{3}}^{2}	\| u_{h} \|_{L_{t}^{\infty} L_{x}^{4}}	+ \| \nabla u_{l} \|_{L_{t,x}^{3}}	\| u_{l} \|_{L_{t}^{3} L_{x}^{12}}^{2}	\| u_{h} \|_{L_{t}^{\infty} L_{x}^{2}}	\lesssim K^{1/3}.
\end{equation}

\noindent Next, by Sobolev embedding,

\begin{equation}\label{4.14}
\int_{J}	\int_{\Omega}	\frac{x_{j}}{(1 + |x|^{2})^{1/2}}	(u_{h}) (\partial_{j} P_{l} F(u))(t,x) dx dt	\lesssim	\| u_{h} \|_{L_{t}^{\infty} L_{x}^{2}}	\| \nabla u_{l} \|_{L_{t,x}^{3}}	\| u_{l} \|_{L_{t}^{3} L_{x}^{12}}^{2}
\end{equation}

\begin{equation}\label{4.15}
+ K^{-1/3}	\| u_{h} \|_{L_{t,x}^{3}}^{2}	\| u_{l} \|_{L_{t}^{3} L_{x}^{12}}	\| u_{l} \|_{L_{t}^{\infty} L_{x}^{4}}	+ K^{-2/3}	\| u_{h} \|_{L_{t,x}^{3}}^{3} \| u_{h} \|_{L_{t}^{\infty} L_{x}^{4}}	\lesssim K^{1/3}.
\end{equation}

\noindent Now consider

\begin{equation}\label{4.16}
\int_{J} \int_{\Omega}	\frac{x_{j}}{(1 + |x|^{2})^{1/2}}	(u_{l}) \partial_{j} (F(u) - F(u_{l}))(t,x) dx dt	+ \int_{J} \int_{\Omega}	\frac{x_{j}}{(1 + |x|^{2})^{1/2}}	P_{l} F(u) \partial_{j} u_{h}(t,x) dx dt.
\end{equation}

\noindent Integrating by parts, since we have already considered $(\ref{4.12})$, $(\ref{4.14})$, it only remains to consider when $\partial_{j}$ hits $\frac{x_{j}}{(1 + |x|^{2})^{1/2}}$. Therefore we have proved

\begin{equation}\label{4.17}
\int_{J}	\int_{\partial \Omega}	|\partial_{n} u(t,x)|^{2} dS_{x} dt		+ \int_{J} \int_{\Omega}	\frac{1}{(1 + |x|^{2})^{3/2}}	[|u_{h}(t,x)|^{2}	+ |\nabla u_{h}(t,x)|^{2}] dx dt
\end{equation}

\begin{equation}\label{4.18}
\int_{J} \int_{\Omega}	\frac{1}{(1 + |x|^{2})^{1/2}}	|u_{h}(t,x)|^{4} dx dt	+ \int_{J} \int_{\Omega}	\frac{1}{(1 + |x|^{2})^{1/2}}	(u_{h}) P_{l} F(u_{h})(t,x) dx dt
\end{equation}

\begin{equation}\label{4.19}
+ \int_{J} \int_{\Omega}	\frac{1}{(1 + |x|^{2})^{1/2}}		O(|u_{h}(t,x)|	|u_{l}(t,x)|^{3})	+ O(|u_{l}(t,x)| |u_{h}(t,x)|^{3}) dx dt	\lesssim K^{1/3}.
\end{equation}

\noindent By Hardy's inequality, that is for $0 \leq s < \frac{d}{2}$,

\begin{equation}\label{4.19.1}
\| \frac{1}{|x|^{s}} f \|_{L^{2}(\mathbf{R}^{d})}	\lesssim_{s,d}	\| |\nabla|^{s} f \|_{L^{2}(\mathbf{R}^{d})},
\end{equation}

\noindent and Sobolev embedding,

\begin{equation}\label{4.20}
\int_{J} \int_{\Omega}	\frac{1}{(1 + |x|^{2})^{1/2}}	(u_{h}) P_{l} F(u_{h})(t,x) dx dt	\lesssim	K^{2/3}	\| \frac{1}{|x|}	u_{h} \|_{L_{t}^{\infty} L_{x}^{2}}	\| u_{h} \|_{L_{t,x}^{3}}^{3}	\lesssim K^{1/3}.
\end{equation}

\noindent Finally,

\begin{equation}\label{4.21}
(\ref{4.19})	\leq	\epsilon	\int_{J} \int_{\Omega}	\frac{1}{(1 + |x|^{2})^{1/2}}	|u_{h}(t,x)|^{4} dx dt	+ C(\epsilon)	\int_{J} \int_{\Omega}	\frac{1}{(1 + |x|^{2})^{1/2}}	|u_{h}(t,x)| |u_{l}(t,x)|^{3} dx dt.
\end{equation}

\begin{equation}\label{4.22}
\int_{J} \int_{\Omega}	\frac{1}{(1 + |x|^{2})^{1/2}}	|u_{h}(t,x)| |u_{l}(t,x)|^{3} dx dt	\lesssim	\| u_{h}(t,x)	\|_{L_{t}^{\infty} L_{x}^{2}}	\| \frac{1}{|x|} |u_{l}|^{3} \|_{L_{t}^{1} L_{x}^{2}}
\end{equation}

\begin{equation}\label{4.23}
\lesssim	\| u_{h}(t,x)	\|_{L_{t}^{\infty} L_{x}^{2}}	\| \nabla u_{l} \|_{L_{t,x}^{3}}	\| u_{l} \|_{L_{t}^{3} L_{x}^{12}}^{2}	\lesssim K^{1/3}.
\end{equation}

\noindent Choosing $\epsilon > 0$ sufficiently small and fixed completes the proof of theorem $\ref{t4.1}$. $\Box$

\begin{theorem}\label{t4.2}
\noindent For $\| u \|_{L_{t,x}^{6}(J \times \Omega)}	= M$, $M$ very large, $u$ solves $(\ref{1.1})$,

\begin{equation}
\| |\nabla|^{-1/2} |u_{h}(t,x)|^{2}	\|_{L_{t,x}^{2}(J \times \Omega)}^{2}	\lesssim	\frac{o(K_{0})}{C_{0}^{3}}.
\end{equation}
\end{theorem}

\noindent $\frac{o(K_{0})}{K_{0}} \rightarrow 0$ as $K_{0} \nearrow \infty$.\vspace{5mm}

\noindent \emph{Proof:} We build on the arguments of \cite{PV}. take the interaction Morawetz quantity

\begin{equation}\label{4.24}
M(t)	= \int_{\Omega \times \Omega}	\frac{(x - y)_{j}}{|x - y|}	Im[\bar{u}_{h}(t,x)	\partial_{j} u_{h}(t,x)] |u_{h}(t,y)|^{2} dx dy.
\end{equation}

\begin{equation}\label{4.25}
\dot{M}(t)	= \int_{\Omega \times \Omega}	\frac{(x - y)_{j}}{|x - y|}	Re[-\partial_{k}^{2} \bar{u}_{h}(t,x) \partial_{j} u_{h}(t,x)	+ \bar{u}_{h}(t,x) \partial_{j} \partial_{k}^{2} u_{h}(t,x)] |u_{h}(t,y)|^{2} dx dy
\end{equation}

\begin{equation}\label{4.26}
+ \int_{\Omega \times \Omega}	\frac{(x - y)_{j}}{|x - y|}	\{ P_{h}(|u|^{2} u), u_{h} \}_{j} (t,x) |u_{h}(t,y)|^{2} dx dy
\end{equation}

\begin{equation}\label{4.27}
- 2 \int_{\Omega \times \Omega}	\frac{(x - y)_{j}}{|x - y|}	Im[\bar{u}_{h}(t,x) \partial_{j} u_{h}(t,x)]	\partial_{k} Im[\bar{u}_{h}(t,y) \partial_{k} u_{h}(t,y)]	dx dy
\end{equation}

\begin{equation}\label{4.28}
+ 2 \int_{\Omega \times \Omega}	\frac{(x - y)_{j}}{|x - y|}	Im[\bar{u}_{h}(t,x) \partial_{j} u_{h}(t,x)]	Im[P_{h}(|u|^{2} u) u_{h}](t,y) dx dy.
\end{equation}

\noindent Integrating by parts, since $u|_{\partial \Omega} = 0$,

\begin{equation}\label{4.29}
(\ref{4.25})	= \int_{\Omega \times \Omega}	\partial_{k}	(\frac{(x - y)_{j}}{|x - y|})	Re[\partial_{k} \bar{u}_{h}(t,x) \partial_{j} u_{h}(t,x)]	|u_{h}(t,y)|^{2} dx dy
\end{equation}

\begin{equation}\label{4.30}
+ \int_{\Omega}		|u_{h}(t,y)|^{2}	\int_{\partial \Omega}	\frac{(x - y)_{j}}{|x - y|}	\nu_{k}	Re[\partial_{k} \bar{u}_{h}(t,x) \partial_{j}(t,x)] dS_{x} dy
\end{equation}

\begin{equation}\label{4.31}
- \int_{\Omega \times \Omega}	\partial_{k} (\frac{(x - y)_{j}}{|x - y|})	Re[\bar{u}_{h}(t,x) \partial_{j} \partial_{k} u_{h}(t,x)]	|u_{h}(t,y)|^{2} dx dy,
\end{equation}

\noindent where $\nu_{k}$ is the outward pointing unit normal to $\partial \Omega$. By theorem $\ref{t2.7}$ $\| u_{h} \|_{L_{t}^{\infty} L_{x}^{2}}	\lesssim	\frac{o(K_{0}^{1/3})}{C_{0}}$, so by theorem $\ref{t4.1}$,

\begin{equation}\label{4.32}
\int_{J} \int_{\Omega}	|u_{h}(t,y)|^{2} \int_{\partial \Omega}	\frac{(x - y)_{j}}{|x - y|}	\nu_{k} Re[\partial_{k} \bar{u}_{h}(t,x) \partial_{j} u_{h}(t,x)] dS_{x} dy	\lesssim \frac{o(K_{0})}{C_{0}^{3}}.
\end{equation}

\noindent Next, integrating by parts,

\begin{equation}\label{4.33}
(\ref{4.31})	= \int_{\Omega \times \Omega}	\partial_{k}(\frac{(x - y)_{j}}{|x - y|})	Re[\partial_{j} \bar{u}_{h}(t,x) \partial_{k} u_{h}(t,x)] |u_{h}(t,y)|^{2} dx dy
\end{equation}

\begin{equation}\label{4.34}
- \frac{1}{2} \int_{\Omega \times \Omega}	(\Delta \Delta |x - y|) |u_{h}(t,x)|^{2} |u_{h}(t,y)|^{2} dx dy.
\end{equation}

\noindent As in $\mathbf{R}^{d}$,

\begin{equation}\label{4.35}
\aligned
\int_{\Omega \times \Omega}	\partial_{k} (\frac{(x - y)_{j}}{|x - y|})	(Re[\partial_{j} \bar{u}_{h}(t,x) \partial_{k} u_{h}(t,x)] |u_{h}(t,y)|^{2}	\\	- Im[\bar{u}_{h}(t,x) \partial_{j} u_{h}(t,x)] Im[\bar{u}_{h}(t,y) \partial_{k} u_{h}(t,y)]) dx dy	\geq 0.
\endaligned
\end{equation}

\noindent Therefore, combining the analysis in theorem $\ref{t4.1}$ with $\| u_{h} \|_{L_{t}^{\infty} L_{x}^{2}}	\lesssim	\frac{o(K_{0}^{1/3})}{C_{0}}$,

\begin{equation}\label{4.36}
\int_{J} \int_{\Omega \times \Omega}	(-\Delta \Delta |x - y|)	|u_{h}(t,x)|^{2} |u_{h}(t,y)|^{2} dx dy dt	+ \int_{J} \int_{\Omega \times \Omega}	\frac{1}{|x - y|} |u_{h}(t,x)|^{4} |u_{h}(t,y)|^{2} dx dy dt
\end{equation}

\begin{equation}\label{4.37}
+ \int_{J} \int_{\Omega \times \Omega}	\frac{(x - y)_{j}}{|x - y|} Im[\bar{u}_{h}(t,x) \partial_{j} u_{h}(t,x)]	Im[P_{h}(|u|^{2} u) \bar{u}_{h}](t,y) dx dy dt	\lesssim	\frac{o(K_{0})}{C_{0}^{3}}.
\end{equation}

\begin{equation}\label{4.38}
Im [P_{h}(|u|^{2} u) \bar{u}_{h}]	= Im[|u|^{4} - |u_{l}|^{4} - (F(u) - F(u_{l})) \bar{u}_{l} - P_{l}(|u|^{2} u) \bar{u}_{h}]	= -Im [(F(u) - F(u_{l})) \bar{u}_{l} + P_{l} F(u) \bar{u}_{h}].
\end{equation}

\noindent By Sobolev embedding

\begin{equation}\label{4.39}
\| u_{h}^{3} u_{l} + P_{l} (u_{h}^{3}) u_{h}	\|_{L_{t,x}^{1}(J \times \Omega)}	\lesssim	\| u_{h} \|_{L_{t,x}^{3}}^{3}	[\| u_{l} \|_{L_{t,x}^{\infty}}	+ K^{-1/3} \| u_{h} \|_{L_{t}^{\infty} L_{x}^{4}}]	\lesssim K^{2/3}.
\end{equation}

\begin{equation}\label{4.40}
\| u_{h}^{2} u_{l}^{2}	+ P_{l}(u_{h}^{2} u_{l}) u_{h} \|_{L_{t,x}^{1}(J \times \Omega)}	\lesssim	\| u_{h} \|_{L_{t,x}^{3}}^{2}	[\| u_{l} \|_{L_{t,x}^{6}}^{2}	+  \| u_{h} \|_{L_{t,x}^{3}}	\| u_{l} \|_{L_{t,x}^{\infty}}]	\lesssim K^{2/3}.
\end{equation}

\begin{equation}\label{4.41}
\| P_{l}(u_{h} u_{l}^{2}) u_{h}	\|_{L_{t,x}^{1}(J \times \Omega)}	\lesssim	\| u_{h} \|_{L_{t,x}^{3}}^{2}	\| u_{l} \|_{L_{t,x}^{6}}^{2}	\lesssim K^{2/3}.
\end{equation}

\noindent Finally, for $u_{h} = \frac{\Delta_{D}}{\Delta_{D}} u_{h}$, integrating by parts,

\begin{equation}\label{4.42}
\int_{J} \int_{\Omega \times \Omega}	\frac{(x - y)_{j}}{|x - y|}	Im[\bar{u}_{h}(t,x) \partial_{j} u_{h}(t,x)]	[u_{h}(t,y) u_{l}(t,y)^{3} + P_{l} F(u_{l})(t,y) u_{h}(t,y)]
\end{equation}

\begin{equation}\label{4.43}
= \int_{J} \int_{\Omega \times \Omega}	\frac{(x - y)_{j}}{|x - y|}	Im[\bar{u}_{h}(t,x) \partial_{j} u_{h}(t,x)]	(\frac{\partial_{k}}{\Delta_{D}}	u_{h}(t,y))	\partial_{k}(u_{l}^{3} + P_{l} F(u_{l}))(t,y) dx dy
\end{equation}

\begin{equation}\label{4.44}
+ \int_{J} \int_{\Omega \times \Omega}	\frac{1}{|x - y|}	Im[\bar{u}_{h}(t,x) \nabla u_{h}(t,x)]	(\frac{\nabla}{\Delta_{D}} u_{h}(t,y))	(u_{l}^{3} + P_{l} F(u_{l}))(t,y) dx dy.
\end{equation}

\begin{equation}\label{4.45}
(\ref{4.43})	\lesssim 	K^{1/3}	\| u_{h} \|_{L_{t}^{\infty} L_{x}^{2}}^{2}	\| \nabla u_{h} \|_{L_{t}^{\infty} L_{x}^{2}}	\| \nabla u_{l} \|_{L_{t,x}^{3}}	\| u_{l} \|_{L_{t}^{3} L_{x}^{12}}^{2}	\lesssim	\frac{o(K_{0})}{C_{0}^{3}}.
\end{equation}

\noindent Now by the Hardy - Littlewood - Sobolev inequality,

\begin{equation}\label{4.46}
(\ref{4.44})	\lesssim	K^{1/3}	\| \nabla u_{h} \|_{L_{t}^{\infty} L_{x}^{2}}	\| u_{h} \|_{L_{t,x}^{3}}^{2}	\| u_{l} \|_{L_{t}^{3} L_{x}^{12}}	\| u_{l} \|_{L_{t}^{\infty} L_{x}^{4}}^{2}	\lesssim \frac{o(K_{0})}{C_{0}^{3}}.
\end{equation}

\noindent This proves theorem $\ref{t4.2}$. $\Box$\vspace{5mm}

\noindent \emph{Proof of theorem $\ref{t4.0}$:} Now we need some constants

\begin{equation}\label{4.47}
0 < \eta_{1} << \eta << 1.
\end{equation}

\noindent Because

\begin{equation}\label{4.48}
K_{0} = C_{0}^{3} \sum_{J_{k} \subset J}	\frac{1}{N_{k}^{3}},
\end{equation}

\noindent for $\eta_{1} > 0$ there exists $K_{0}(\eta_{1})$ sufficiently large such that there exists $J_{k} \subset J$ with

\begin{equation}\label{4.49}
\| |\nabla|^{-1/2} |u_{h}|^{2}	\|_{L_{t,x}^{2}(J_{k} \times \Omega)}^{2}	\leq \frac{\eta_{1}}{N_{k}^{3}}.
\end{equation}

\noindent Now on each $J_{k}$,

\begin{equation}\label{4.50}
\| \nabla |u_{h}|^{2}	\|_{L_{t,x}^{2}(J_{k} \times \Omega)}	\lesssim	\| \nabla u_{h} \|_{L_{t,x}^{3}}	\| u_{h} \|_{L_{t,x}^{6}}	\lesssim 1.
\end{equation}

\noindent Therefore, by interpolation

\begin{equation}\label{4.51}
\| u_{h} \|_{L_{t,x}^{4}(J_{k} \times \Omega)}^{4}	\lesssim	\frac{\eta_{1}^{2/3}}{N_{k}^{2}}.
\end{equation}

\noindent Now by interpolation, since $\| u \|_{L_{t}^{2+} L_{x}^{\infty-}(J_{k} \times \Omega)}	\lesssim 1$, $\| u \|_{L_{t,x}^{6}(J_{k} \times \Omega)}	\gtrsim 1$,

\begin{equation}\label{4.52}
\| u \|_{L_{t}^{\infty} L_{x}^{4}(J_{k} \times \Omega)}	\gtrsim 1.
\end{equation}

\noindent Therefore there exists $t_{k} \in J_{k}$ such that $\| u(t_{k}) \|_{L_{x}^{4}(\Omega)} \gtrsim 1$. Moreover, by theorem $\ref{t2.7}$ and Sobolev embedding,

\begin{equation}\label{4.53}
\| u_{\frac{1}{C(\eta)} N_{k} \leq \cdot \leq C(\eta) N_{k}}(t_{k}) \|_{L_{x}^{4}(\Omega)}	\gtrsim 1.
\end{equation}

\noindent Take $K_{0}$ sufficiently large so that $K_{0}^{-1/3} << \frac{1}{C(\eta)}$.

\begin{equation}\label{4.54}
\frac{d}{dt} \int_{\Omega}	|u_{\frac{1}{C(\eta)} N_{k} \leq \cdot \leq C(\eta) N_{k}}(t,x)|^{4} dx	\lesssim C(\eta)^{2} N_{k}^{2}.
\end{equation}

\noindent Therefore, for some $\delta > 0$, $\| u_{h}(t) \|_{L_{x}^{4}(\Omega)} \gtrsim 1$ on $[t_{k} - \frac{\delta}{C(\eta)^{2} N_{k}^{2}}, t_{k} + \frac{\delta}{C(\eta)^{2} N_{k}^{2}}]$. However this implies

\begin{equation}\label{4.55}
\| u_{h} \|_{L_{t,x}^{4}(J_{k} \times \Omega)}^{4} \gtrsim \frac{1}{C(\eta)^{2} N_{k}^{2}},
\end{equation}

\noindent which contradicts $(\ref{4.51})$. $\Box$\vspace{5mm}

\noindent Theorem $\ref{t3.0}$ combined with theorem $\ref{t4.0}$ proves theorem $\ref{t2.8}$. It only remains to prove $(\ref{4.5})$.

\section{Endpoint argument}
\noindent It only remains to prove $(\ref{4.5})$. To do this we will upgrade lemma $\ref{l3.2}$ to involve $l^{2}$ summation. For $K = \sum_{J_{k} \subset J} \frac{1}{N_{k}^{3}}$ we define the norm

\begin{equation}\label{5.1}
\| u \|_{X(J \times \Omega)}^{2}	\equiv	\sum_{K^{-1/3} \leq N_{j}}	\frac{1}{K N_{j}}		\sum_{J_{N_{j}}^{l} \subset J}	\| u_{N_{j}}	\|_{U_{\Delta}^{2}(J_{N_{j}}^{l} \times \Omega)}^{2}.
\end{equation}

\begin{theorem}\label{t5.1}
If $u$ is a solution to $(\ref{1.1})$, $\| u \|_{L_{t,x}^{6}(J \times \Omega)}	= M$ for some $M$ sufficiently large and fixed, $E(u(t)) = E_{0}$, then

\begin{equation}\label{5.2}
\| u \|_{X(J \times \Omega)}	\lesssim 1.
\end{equation}
\end{theorem}

\noindent \emph{Proof:} We again take $(\ref{3.7})$. First consider the bad intervals $J_{N_{j},b}^{l}$. By lemma $\ref{l3.2}$,

\begin{equation}\label{5.3}
 \sum_{K^{-1/3} \leq N_{j}}	\frac{1}{K N_{j}}	\sum_{J_{N_{j},b}^{l}}	\| P_{N_{j}} u \|_{U_{\Delta}^{2}(J_{N_{j},b}^{l} \times \Omega)}^{2}
 \lesssim \frac{1}{K} \sum_{J_{k} \subset J}	\frac{1}{c^{3} N_{k}^{3}}	\lesssim	\frac{1}{c^{3}}.
\end{equation}

\noindent Now turn to the good intervals.

\begin{equation}\label{5.5}
\| \Delta P_{N_{j}} F(u_{\leq N_{j}}) \|_{DU_{\Delta}^{2}(J_{N_{j}, g}^{l} \times \Omega)}		\lesssim	\| \Delta u_{\leq N_{j}}	\|_{L_{t}^{2+} L_{x}^{4-}(J_{N_{j}, g}^{l} \times \Omega)}	\| \nabla u_{\leq N_{j}} \|_{L_{t}^{2+} L_{x}^{4-}}	\| u_{\leq N_{j}}	\|_{L_{t}^{\infty} \dot{H}_{x}^{1}}
\end{equation}

\begin{equation}\label{5.6}
\lesssim	\eta \| \Delta u_{\leq N_{j}} \|_{U_{\Delta}^{2}(J_{N_{j}, g}^{l} \times \Omega)}.
\end{equation}

\noindent The last inequality follows from lemma $\ref{l3.2}$. Now we use the fact that an interval $J_{N}^{l}$, $N \leq N_{j}$ overlaps $\lesssim (\frac{N_{j}}{N})^{3}$ intervals $J_{N_{j},g}^{l}$.

\begin{equation}\label{5.7}
\frac{1}{K}	\sum_{K^{-1/3} \leq N_{j}}	\frac{1}{N_{j}}	\sum_{J_{N_{j},g}^{l} \subset J}	\| P_{N_{j}} F(u_{\leq N_{j}}) \|_{DU_{\Delta}^{2}(J_{N_{j},g}^{l} \times \Omega)}^{2}
\end{equation}

\begin{equation}\label{5.8}
\lesssim	\frac{\eta}{K}	\sum_{K^{-1/3} \leq N_{j}}		\sum_{K^{-1/3} \leq N \leq N_{j}}	\frac{N^{4}}{N_{j}^{5}}	(\frac{N_{j}}{N})^{3}	\sum_{J_{N}^{l} \subset J}	\| u_{N} \|_{U_{\Delta}^{2}(J_{N}^{l} \times \Omega)}^{2}
\end{equation}

\begin{equation}\label{5.9}
+ \eta	\sum_{K^{-1/3} \leq N_{j}}	\frac{K^{-2/3}}{N_{j}^{2}}	\| \nabla u_{\leq K^{-1/3}} \|_{U_{\Delta}^{2}(J \times \Omega)}^{2}	\lesssim \eta (1 + \| u \|_{X(J \times \Omega)}^{2}).
\end{equation}

\noindent Next,

\begin{equation}\label{5.10}
\| P_{N_{j}} F(u_{N_{j} < \cdot c N_{k}}) \|_{DU_{\Delta}^{2}(J_{N_{j}}^{l} \times \Omega)}	\lesssim	\| P_{N_{j}} F(u_{N_{j} < \cdot < cN_{k}}) \|_{L_{t}^{2-} L_{x}^{\frac{4}{3}+}(J_{N_{j}, g}^{l} \times \Omega)}
\end{equation}

\begin{equation}\label{5.11}
\lesssim	N_{j}	\sum_{N_{j} < \cdot M < c N_{k}}	\| u_{M} \|_{L_{t}^{2+} L_{x}^{4-}(J_{N_{j},g}^{l} \times \Omega)}^{1+}	\| u_{M} \|_{L_{t}^{\infty} L_{x}^{2}(J_{N_{j}, g}^{l} \times \Omega)}^{1-} \| u_{M} \|_{L_{t}^{\infty} L_{x}^{4}(J_{N_{j},g}^{l} \times \Omega)}
\end{equation}

\begin{equation}\label{5.12}
\lesssim	\sum_{N_{j} < M < cN_{k}} \eta (\frac{N_{j}}{M})^{1-} \| u_{M} \|_{U_{\Delta}^{2}(J_{N_{j}}^{l} \times \Omega)}.
\end{equation}

\noindent The last inequality follows from lemma $\ref{l3.2}$. Now by taking the convolution of an $L^{1}$ function with an $L^{2}$ function,

\begin{equation}\label{5.13}
\frac{\eta^{1-}}{K}	\sum_{K^{-1/3} \leq N_{j}}	\sum_{J_{N_{j},g}^{l}		\subset J}	(\sum_{N_{j} < M < cN_{k}}	\frac{1}{N_{j}^{1/2}}	(\frac{N_{j}}{M})^{1-} (\sum_{J_{M}^{l'} \cap J_{N_{j}, g}^{l} \neq \emptyset}	\| u_{M} \|_{U_{\Delta}^{2}(J_{M}^{l'} \times \Omega)}^{2})^{1/2})^{2}
\end{equation}

\begin{equation}\label{5.14}
\lesssim	\frac{\eta^{1-}}{K}	\sum_{K^{-1/3} \leq N_{j}}	\sum_{J_{N_{j},g}^{l} \subset J}		(\sum_{N_{j} < M < cN_{k}}	(\frac{N_{j}}{M})^{\frac{1}{2}-}	(\sum_{J_{M}^{l'} \cap J_{N_{j},g}^{l}} \frac{1}{M}	\| u_{M} \|_{U_{\Delta}^{2}(J_{M}^{l'} \times \Omega)}^{2})^{1/2})^{2}
\end{equation}

\begin{equation}\label{5.15}
\lesssim	\frac{\eta^{1-}}{K}	\sum_{K^{-1/3} \leq M}			(\sum_{J_{M}^{l'} \subset J}	\frac{1}{M} \| u_{M} \|_{U_{\Delta}^{2}(J_{M}^{l'} \times \Omega)}^{2})	\lesssim	\eta^{1-}	\| u \|_{X(J \times \Omega)}^{2}.
\end{equation}

\noindent Finally,

\begin{equation}\label{5.16}
\| \nabla P_{N_{j}} F(u_{> cN_{k}}) \|_{DU_{\Delta}^{2}(J_{N_{j},g}^{l} \times \Omega)}	\lesssim	\| \nabla P_{N_{j}} F(u_{> cN_{k}}) \|_{L_{t}^{2-} L_{x}^{\frac{4}{3}+}(J_{N_{j},g}^{l} \times \Omega)}.
\end{equation}

\noindent By lemma $\ref{l3.2}$,

\begin{equation}\label{5.17}
\| \nabla P_{N_{j}} F(u_{> cN_{k}}) \|_{L_{t}^{2-} L_{x}^{\frac{4}{3}+}(J_{N_{j},g}^{l} \times \Omega)}	\lesssim 1.
\end{equation}

\noindent Therefore

\begin{equation}\label{5.18}
\sum_{K^{-1/3} \leq N_{j}}	\frac{1}{K N_{j}^{3}}	\sum_{J_{N_{j},g}^{l} \subset J}	\| \nabla P_{N_{j}} F(u_{> cN_{k}}) \|_{L_{t}^{2-} L_{x}^{\frac{4}{3}+}(J_{N_{j},g} \times \Omega)}^{2}
\end{equation}

\begin{equation}\label{5.19}
\lesssim	\sum_{K^{-1/3} \leq N_{j}}	\frac{1}{K N_{j}^{3}}	\sum_{J_{k} \subset J : 2 N_{j} \leq N_{k}}	(N_{j}^{2}	\frac{1}{(cN_{k})^{2})^{2-}}
\end{equation}

\begin{equation}\label{5.20}
\lesssim	\frac{1}{K}	\sum_{J_{k}}	\frac{1}{(cN_{k})^{4-}}		\sum_{N_{j} \leq \frac{N_{k}}{2}}	N_{j}^{1-}	\lesssim	\frac{1}{c^{4-}}.
\end{equation}

\noindent Finally consider the term

\begin{equation}\label{5.21}
\| P_{N_{j}} ((u_{< N_{j}})(u_{> N_{j}}^{2})	+ (u_{< N_{j}}^{2})(u_{> N_{j}})) \|_{DU_{\Delta}^{2}(J_{N_{j},g}^{l} \times \Omega)}.
\end{equation}

\noindent First take $(u_{< N_{j}}^{2}) u_{> N_{j}}$. Suppose $v = [\Psi(\frac{1}{4} N_{j}^{-2} \Delta_{D}) + \Psi(N_{j}^{-2} \Delta_{D}) + \Psi(4 N_{j}^{-2} \Delta_{D})]v$, $\| v \|_{V_{\Delta}^{2}(J_{N_{j},g}^{l} \times \Omega} = 1$.

\begin{equation}\label{5.22}
\| (u_{N}) v (u_{< N_{j}})(u_{> N_{j}})	\|_{L_{t,x}^{1}(J_{N_{j},g}^{l} \times \Omega)}	\lesssim	\| v (u_{N}) \|_{L_{t,x}^{2}}	\| u_{< N_{j}} \|_{L_{t}^{2+} L_{x}^{\infty-}}	\| u_{> N_{j}} \|_{L_{t}^{\infty-} L_{x}^{2+}}.
\end{equation}

\noindent By lemma $\ref{l3.2}$, theorem $\ref{t3.1}$,

\begin{equation}\label{5.23}
\| u_{> N_{j}} \|_{L_{t}^{4} L_{x}^{8/3}}	\lesssim	\frac{1}{N_{j}},
\end{equation}

\begin{equation}\label{5.24}
 \| u_{< N_{j}} \|_{L_{t}^{4} L_{x}^{8}}	\lesssim \eta^{1/2}.
\end{equation}

\begin{equation}\label{5.25}
 \| (u_{N}) v \|_{L_{t,x}^{2}}	\lesssim	\| (u_{N})v \|_{L_{t}^{2+} L_{x}^{2-}}^{1-}	\| u_{N} \|_{L_{t}^{2+} L_{x}^{\infty-}}^{+}	\| v \|_{L_{t}^{2+} L_{x}^{4-}}^{+}	\lesssim	\frac{N^{3/2-}}{N_{j}^{1/2-}}	\| u_{N} \|_{U_{\Delta}^{2}(J_{N_{j},g}^{l} \times \Omega)}.
\end{equation}

\noindent The last inequality follows from interpolating theorem $\ref{t1.2}$ with theorem$\ref{t2.7}$, $V_{\Delta}^{2} \subset U_{\Delta}^{p}$ for $p > 2$.

\begin{equation}\label{5.26}
\frac{\eta^{1/2}}{K}	\sum_{K^{-1/3} \leq N_{j}} \frac{1}{N_{j}}  \sum_{J_{N_{j},g}^{l} \subset J} \sum_{K^{-1/3} \leq N \leq N_{j}} (\frac{N}{N_{j}})^{3-} \| u_{N} \|_{U_{\Delta}^{2}(J_{N_{j},g}^{l} \times \Omega)}^{2}	\lesssim	\eta^{1/2} \| u \|_{X(J \times \Omega)}^{2}.
\end{equation}

\noindent By lemma $\ref{l3.2}$,

\begin{equation}\label{5.27}
\frac{\eta^{1/2}}{K}	\sum_{K^{-1/3} \leq N_{j}} \frac{1}{N_{j}} K N_{j}^{3} \frac{K^{-1/3}}{N_{j}^{3}} \| \nabla u_{K^{-1/3}} \|_{U_{\Delta}^{2}(J \times \Omega)}^{2} \lesssim 1.
\end{equation}

\noindent Finally,

\begin{equation}\label{5.28}
\| (u_{> N_{j}})^{2} (u_{< N_{j}}) \|_{L_{t}^{2-} L_{x}^{4/3+}(J_{N_{j},g}^{l} \times \Omega)}	\lesssim	N_{j} \| u_{< N_{j}} \|_{L_{t}^{\infty-} L_{x}^{4+}} \sum_{M \geq N_{j}} \| u_{M} \|_{L_{t}^{2+} L_{x}^{4-}} \| u_{M} \|_{L_{t}^{\infty} L_{x}^{2}}
\end{equation}

\begin{equation}\label{5.29}
 \lesssim	\eta^{1-} \sum_{M \geq N_{j}} \frac{N_{j}}{M} \| u_{M} \|_{U_{\Delta}^{2}(J_{N_{j},g}^{l} \times \Omega)}.
\end{equation}

\begin{equation}\label{5.30}
 \frac{\eta^{1-}}{K} \sum_{K^{-1/3} \leq N_{j}} \frac{1}{N_{j}} \sum_{J_{N_{j},g}^{l} \subset J} (\sum_{M \geq N_{j}} \frac{N_{j}}{M} \| u_{M} \|_{U_{\Delta}^{2}(J_{N_{j},g}^{l} \times \Omega)})^{2} \lesssim \eta^{1-} \| u \|_{X(J \times \Omega)}^{2}.
\end{equation}

\noindent Finally, for $(\ref{3.7})$ choose $t_{N_{j}}^{l} \in J_{N_{j},g}^{l}$ such that

\begin{equation}\label{5.31}
\| \nabla P_{N_{j}} u(t_{N_{j}}^{l}) \|_{L_{x}^{2}(\Omega)}	= \inf_{t \in J_{N_{j},g}^{l}}	(\| \nabla P_{N_{j}} u(t) \|_{L_{x}^{2}(\Omega)}).
\end{equation}

\noindent This implies

\begin{equation}\label{5.32}
\sum_{K^{-1/3} \leq N_{j}}	\frac{1}{K N_{j}}		\sum_{J_{N_{j},g}^{l} \subset J}		\| P_{N_{j}} u(t_{N_{j}}^{l}) \|_{L_{x}^{2}(\Omega)}^{2}	\lesssim	1.
\end{equation}

\noindent Therefore, we have proved

\begin{equation}\label{5.33}
 \| u \|_{X(J \times \Omega)}^{2}	\lesssim 1 + \eta^{1/2} \| u \|_{X(J \times \Omega)}^{2}.
\end{equation}

\noindent This completes the proof of theorem $\ref{t5.1}$. $\Box$\vspace{5mm}

\noindent By lemma $\ref{l3.2}$,

\begin{equation}\label{5.35}
\| u_{N} \|_{L_{t}^{2+} L_{x}^{4-}(J \times \Omega)}	\lesssim	 (\sum_{J_{N}^{l} \subset J}	\| u_{N} \|_{L_{t}^{2+} L_{x}^{4-}(J_{N}^{l} \times \Omega)}^{2+})^{1/2+}	\lesssim	N^{-\frac{1}{2} + \frac{1}{2+}}(\sum_{J_{N}^{l} \subset J}	\| u_{N} \|_{U_{\Delta}^{2}(J_{N}^{l} \times \Omega)}^{2})^{1/(2+)}.
\end{equation}

\noindent Also,

\begin{equation}\label{5.36}
\| u_{N} \|_{L_{t}^{\infty-} L_{x}^{2+}(J \times \Omega)}	\lesssim	M^{\frac{1}{\infty-} - 1}  (\sum_{J_{N}^{l} \subset J}	\| u_{N} \|_{U_{\Delta}^{2}(J_{N}^{l} \times \Omega)}^{2})^{1/(\infty-)}
\end{equation}

\noindent Therefore,

\begin{equation}\label{5.37}
 \sum_{K^{-1/3} \leq N_{1} \leq N_{2} \leq N_{3}}	\| u_{N_{1}} \|_{L_{t}^{2+} L_{x}^{4-}(J \times \Omega)}	\| u_{N_{2}} \|_{L_{t}^{2+} L_{x}^{4-}(J \times \Omega)} \| u_{N_{3}} \|_{L_{t}^{\infty-} L_{x}^{2+}(J \times \Omega)}
\end{equation}

\begin{equation}\label{5.38}
\aligned
\lesssim	 \sum_{K^{-1/3} \leq N_{1} \leq N_{2} \leq N_{3}}	N_{1}^{-1/2 + \frac{2}{2+}}		N_{2}^{-1/2 + \frac{2}{2+}}		N_{3}^{\frac{2}{\infty-}	- 1}	 (\sum_{J_{N_{3}}^{l} \subset J}	\| u_{N_{3}} \|_{U_{\Delta}^{2}(J_{N_{3}}^{l} \times \Omega)}^{2})^{1/(\infty-)}	\\
\times	(\sum_{J_{N_{1}}^{l} \subset J}	\| u_{N_{1}} \|_{U_{\Delta}^{2}(J_{N_{1}}^{l} \times \Omega)}^{2})^{1/(2+)}	(\sum_{J_{N_{2}}^{l} \subset J}	\| u_{N_{2}} \|_{U_{\Delta}^{2}(J_{N_{2}}^{l} \times \Omega)}^{2})^{1/(2+)}	\lesssim	\| u \|_{X(J \times \Omega)}^{2}.
\endaligned
\end{equation}

\noindent The proof of theorem $\ref{t1.1}$ is now complete.

\newpage

\nocite*
\bibliographystyle{plain}
\bibliography{energy}

\end{document}